\documentclass[12pt]{article}

\usepackage{graphicx,amssymb,amsmath,amsthm,amsfonts}
\usepackage{subfigure,epsfig}
\usepackage{setspace}
\usepackage{natbib}
\usepackage{authblk,color}
\usepackage[outdir=./]{epstopdf}

\theoremstyle{definition}

\title{\bfseries A finite mixture model approach to regression under covariate misclassification}

\author[1]{P. Richard Hahn\thanks{P. Richard Hahn, University of Chicago Booth School of Business,
5807 S. Woodlawn Avenue, Chicago, IL 60637, U.S.A. Email: richard.hahn@chicagobooth.edu}}
\author[2]{Michelle Xia}
\affil[1]{Booth School of Business, University of Chicago, Chicago, IL, 60637, U.S.A.}
\affil[2]{Division of Statistics, Northern Illinois University, Dekalb, IL, 60115, U.S.A.} 

\begin{document}
\DeclareGraphicsExtensions{.eps} 

\maketitle
  
\subparagraph{Abstract.} This paper considers the problem of mismeasured categorical covariates in the context of regression modeling; if unaccounted for, such misclassification is known to result in misestimation of model parameters.  Here, we exploit the fact that explicitly modeling covariate misclassification leads to a mixture representation. Assuming common parametric families for the mixture components, and assuming that the misclassification occurrence is independent of the response variable, the mixture representation permits model parameters to be identified even when misclassification probabilities are unknown. Previous approaches to covariate misclassification use multiple surrogate covariates and/or validation data on the magnitude of errors. Based on this mixture structure, we demonstrate that valid inference can be performed on all the parameters even when no such additional information is available. Using Bayesian inference, the method allows for learning from data combined with external information on the magnitude of errors when such information does become available. The method is applied to adjust for misclassification on self-reported cocaine use in the Longitudinal Studies of HIV-Associated Lung Infections and Complications (Lung HIV). We find a substantial and statistically significant effect of cocaine use on pulmonary complications measured by the relative area of emphysema, whereas a regression that does not adjust for misclassification yields a much smaller estimate.



\section{Introduction}
{\em Misclassification} refers to measurement error on categorical or binary variables, so that some observations within a data set might record a category other than the truth. For example, our motivating example arises because survey respondents are often reluctant to truthfully report drug use. When covariates or exposure variables are subject to misclassification, naive (unadjusted) estimation of regression effect can be biased (see, e.g., \cite{liu1991, beavers2012bayesian, wang14}). 

Previous methods that adjust regression estimates for covariate misclassification {use} information from either multiple surrogates of the misclassified variable (see, e.g., \cite{liu1991}) or validation data on the misclassification probabilities (see, e.g., \citet*{Chu10}). Recently, \cite{Shieh09} and \cite{hubbardaccounting} studied a mixture representation for regression with misclassified categorical covariates, respectively for normal and binary response variables. Based on the mixture representation, the papers proposed frequentist methods adjusting for the misclassification, for situations when the reclassification probabilities are known (e.g., from previous studies).

In this paper, we expand the methods of \cite{Shieh09} and \cite{hubbardaccounting} to study the general structure of mixture regression models resulting from covariate misclassification. The mixture representation demonstrates that model parameters can be identified even when none of the aforementioned side information is available. Specifically, the proposed method does not rely on additional information on the magnitude of errors, such as the reclassification probabilities assumed in \cite{Shieh09} and \cite{hubbardaccounting}. Using Bayesian inference, external information on the reclassification probabilities can be incorporated through priors, when such information does become available. This helps us benefit from the capacity of existing methods in incorporating external information, while allowing additional statistical learning from data. Our approach can be applied to a wide class of regression models, in addition to the normal and binary cases studied in earlier papers, including generalized linear models and parametric survival models, provided the mixture components are identifiable. Our model allows for additional (correctly measured) covariates, provided that the misclassification is {\em nondifferential}, meaning that misclassification is conditionally independent of the response variable. We study the effectiveness of the mixture representation approach via an asymptotic analysis of efficiency and via a Monte Carlo study. The upshot of these investigations is that, if the effect size is large enough, our approach gives up little statistical efficiency relative to the case where the classification probabilities are known {\em a prior}. When the effect size is small, we demonstrate that the proposed model works as good as alternatives in cases where accurate information is available on the magnitude of errors. When the information is inaccurate, our method is able to learn this from the data, thereby reducing bias in the estimation. 


Our motivating empirical example is the multi-center Lung HIV study (\cite{crothers2011hiv}) that was conducted during 2007 to 2013 by the National Heart, Lung, and Blood Institute (NHLBI). The goal of this study was to understand the relationship between human immunodeficiency virus (HIV) infection and pulmonary diseases. Data from the Lung HIV study and its sub-studies have been analyzed in the medical literature (\cite{drummond2015factors, leader2016risk}) in an effort to identify risk factors for pulmonary complications related to HIV infection. Other papers (e.g., \cite{sigel2014findings,depp2016risk}) assess HIV infection itself as a risk factor for pulmonary condition; \cite{drummond2013effect,lambert2015hiv} focus on injection drug users. 

For the Lung HIV study, self-reported traits were collected on use of illegal drugs such as crystal methamphetamine, crack cocaine, marijuana and heroin. Previous literature indicated a strong clinical connection between pulmonary diseases and use of recreational drugs, particularly cocaine use (\citet*{yakel1995pulmonary,alnas2010clinical}; \cite{fiorelli2016spontaneous}). On the relationship between drug use and pulmonary function specifically in HIV-infected individuals,  \cite{simonetti2014pulmonary} reported no significant effect of recreational drug use on pulmonary function, based on two sub-studies of the Lung HIV data (Multicenter AIDS Cohort Study and Women's Interagency HIV Study). The goal of our analysis is to revisit this question using the full Lung HIV study, while accounting for the covariate misclassification due to the probable inaccuracy of self-reported drug use. Due to the anticipated large effect size of cocaine use and the unavailability of external information, our method becomes a natural choice in adjusting for misclassification when assessing the cocaine effect.

The rest of the paper is organized as follows. In Section \ref{SocSec0}, we provide a general framework for regression models adjusting for covariate misclassification. In Section \ref{SecEst}, we obtain and evaluate the asymptotic results on the efficiency loss for the normal model. In Section \ref{Simu}, we present a simulation study on finite samples. In Section \ref{appl}, the methodology is applied to adjustment for misclassification in the self-reported cocaine use, when assessing its effect on lung density measures. Section \ref{concl} concludes the paper.

\section{A finite mixture model representation}\label{SocSec0}

Throughout the paper, we use upper-case Roman letters to denote random variables such as $Y$, with the corresponding lower-case letter $y$ denoting an observed value of the variable. We use lower-case bold letters such as $\mathbf{x}$ to denote a vector of random variables, and the notation $\boldsymbol{x}$ for an observed value of the random vector. Similarly, the bold version of parameters such as $\boldsymbol{\alpha}$ represents a vector of parameters. All vectors are column vectors, unless specified otherwise. Upper-case bold letters such as $\mathbf{P}$ are used to denote a matrix. For notational simplicity, we present a population model based on a single set of random variables, without involving an observed sample of size $n>1$.

\subsection{A finite mixture regression representation}
Let $V$ be a categorical variable taking $K\geq 2$ categories, which we will denote (0, 1, $\cdots$, $K-1$), and let $P(V=k)=\pi_k \geq 0$ with $\sum_k\pi_k=1$ denote the corresponding category probabilities. Denote by $V^*$ the observed version of $V$ with misclassification. The \emph{misclassification probability} is defined as $p_{kj}=\mbox{P}(V^*=j\,|\,V=k)$, with $\sum_jp_{kj}=1$ for $k=0,\,1,\,\cdots,\,K-1$. Hence, the severity of misclassification can be represented with the classification matrix $\mathbf{P}=(p_{kj})$. Interested readers may refer to \cite{buonaccorsi10,GustGr15,Yi2016} for comprehensive reviews on the issue and adjustment of misclassification in categorical variables.

Here, we formulate the misclassification problem in the context of parametric regression, where the distribution of the response variable $Y$ is modeled conditional on a categorical covariate $V$ and additional (accurately measured) covariates denoted using a vector $\mathbf{x}$. We write the probability function of $(Y\,|\,V,\,\mathbf{x})$ as $f_{Y}(y\,|\,\boldsymbol{\alpha},\,\boldsymbol{\beta},\,\boldsymbol{\varphi},\,V,\,\mathbf{x})$ for parameters $(\boldsymbol{\alpha},\,\boldsymbol{\beta},\,\boldsymbol{\varphi})$. Specifically, in the case of generalized linear models (GLMs) with a binary $V$, we have $g(\mbox{E}(Y\,|\,V,\,\mathbf{x}))=\alpha_0+\alpha_1V + \mathbf{x}'\boldsymbol{\beta}$, with $\boldsymbol{\alpha}=(\alpha_0, \,\alpha_1)'$, a vector $\boldsymbol{\beta}$ containing the regression coefficients for $\mathbf{x}$, $\boldsymbol{\varphi}$ containing possible nuisance parameters such as a dispersion parameter, and $g(\cdot)$ being the link function.

Assume that the misclassification is nondifferential on $Y$, meaning that the occurrence of misclassification is conditionally independent of $Y$. Given the observed covariates $V^{*}$ and $\mathbf{x}$, the conditional distribution of $Y$ has mixture representation
\begin{align}\label{mix.eq1}
 f_Y(y\,|\,V^{*},\,\mathbf{x})&=\sum_{j=0}^{K-1}q_{j}(V^{*},\,\mathbf{x})\,f_Y(y\,|\,\boldsymbol{\alpha},\,\boldsymbol{\beta},\,\boldsymbol{\varphi},\,V=j,\,\mathbf{x}),
\end{align} 
where the mixture weights are the {\em reclassification probabilities}, defined as $q_{kj}(V^*,\mathbf{x})=\mbox{P}(V=j\,|\,V^*=k,\,\mathbf{x})$; $k=0,\,1,\,\cdots,\,K-1$. The specific form of the $q_{j}(\cdot)$, which is a function from the domain of $(V^{*},\,\mathbf{x})$ to $[0,1]$, depends on the joint distribution of $(V,\,V^{*},\,\mathbf{x})$. The $K$-by-$K$ {\em reclassification matrix} is denoted $\mathbf{Q}=(q_{kj})$, with $\sum_jq_{kj}=1$ for $k=0,\,1,\,\cdots,\,K-1$, where dependence on $(V^{*},\,\mathbf{x})$ is left implicit. Equation (\ref{mix.eq1}) gives a weighted sum of $K$ regression models, each corresponding to a particular unobserved value of $V$, with corresponding weights depending on the observed value of $V^*$ (which can take the same $K$ distinct values as $V$). This extends the mixture distribution representation in \cite{Shieh09} to a more general \emph{mixture regression model} framework (see, e.g., \cite{grun2008}).


\subsection{Identification}

The model in (\ref{mix.eq1}) is a special case of a more general class of models known as mixture regression models with concomitant variables (\cite{grun2008}) or mixture-of-experts models (\cite{jacobs1991adaptive}). This more general representation takes the form 

\begin{align}\label{mix.eq2}
 f_Y(y\,|\,\mathbf{w},\,\mathbf{x})&=\sum_{j=0}^{K-1}\phi_{j}(\mathbf{w})\,f_Y(y\,|\,\boldsymbol{\omega}_j,\,\mathbf{w},\,\mathbf{x}),
\end{align} 
where each $\boldsymbol{\omega}_j$ is a vector of component-specific parameters, and $\mathbf{w}$ and $\mathbf{x}$ are vectors of observed covariates. According to the first nomenclature, covariates $\mathbf{w}$ appearing in the mixture weights $\phi_{j}(\cdot)$ are called concomitant variables; according to the second, each regression component is referred to as an ``expert'' and the $\phi_{j}(\cdot)$ are referred to as gating functions. Although it is possible that the gating functions take more general forms, it is common to use a multinomial logit form
\begin{align}\label{mix.eq3}
\phi_{j}(\mathbf{w})&= \frac{\exp{(\nu_j + \mathbf{w}'\boldsymbol{\gamma}_j)}}{\sum_{h=0}^{K-1} \exp{(\nu_h + \mathbf{w}'\boldsymbol{\gamma}_h)}},\quad j=0,\,1,\,\cdots,\,K-1.
\end{align} 

Identification of parameters in these models is a nontrivial issue, which is discussed in \cite{hennig2000identifiability}, \cite{grun2008finite}, and \cite{jiang1999identifiability}. In these papers it is shown that models of the form in (\ref{mix.eq2}) are identifiable --- meaning there is a unique mapping from probability functions  $f_Y(y\,|\,V^{*},\,\mathbf{x})$ to parameters $\boldsymbol{\Omega} = (\boldsymbol{\omega}_0, \boldsymbol{\omega}_1, \dots, \boldsymbol{\omega}_{K-1})$ --- provided that several criteria are satisfied. We discuss these conditions next, explaining why the models studied in this paper satisfy them.

First, the family of component densities $f_Y(\cdot)$ must be identifiable. That is, given a finite mixture of distributions in this family, one must be able to uniquely decompose it into its constituent components. This is known to be true of many common densities, including the normal, gamma, and Poisson distributions (\citet*{yakowitz1968,atienza2006}; \cite{atienza2007}) considered in this paper. Note that these results require that the mixture representation is {\em irreducible}, meaning that all $K$ mixture components are distinct; this is true for our model if and only if $\alpha_1 \neq 0$.

Second, one must be able to order the parameter vectors $(\boldsymbol{\omega}_0 \prec \boldsymbol{\omega}_1  \prec \dots, \boldsymbol{\omega}_{K-1})$ uniquely. This ordering serves to (arbitrarily) break the symmetry inherent to mixture representations, which is that they are invariant under permutations of the component labels. \cite{jiang1999identifiability} provides a recipe for establishing such an ordering for exponential family models. Moreover, this condition is easily satisfied for (\ref{mix.eq1}) because only the (scalar) $\alpha_1 V$ term varies by component. 

Regarding this second condition, it is important to note that in the misclassified covariate context the parameter ordering is not arbitrary. For example, in the case of a binary variable $V$, one component of the mixture model corresponds to the mean response when $V = 0$ and the other to the mean response when $V = 1$; clearly these designations have applied meaning and are not interchangeable. Fortunately, the relative ordering of these components will be preserved as long as the misclassification is not {\em systematic}, meaning that the probability of misclassification is higher than that of correct classification and that the direction (sign) of the slope parameter $\alpha_1$ is known. See  \citet*{weinberg1994} for additional discussion.

Third, one must parametrize the reclassification probabilities in terms of a {\em base category}, for example by setting $\nu_{K-1} = 0$ and $\boldsymbol{\gamma}_{K-1} = \mathbf{0}$. This breaks the invariance to location shifts of the gating parameters $(\nu_j, \boldsymbol{\gamma}'_j)$.

Finally, \cite{hennig2000identifiability} points out that a mixture of regression model has additional criteria that must be satisfied for a given set of covariate values. There and in \cite{grun2008finite}, concrete examples are produced which show that different parameters $\boldsymbol{\Omega} \neq \boldsymbol{\Omega}'$ can give the same likelihood evaluation if ``[component] labels are fixed in one covariate point according to some ordering constraint, [but] labels switch in other covariate points for different parameterizations of the model.'' Note that \cite{jiang1999identifiability} rule out this possibility by stipulating that $\mathbf{w} = \mathbf{x}$ defines a set with a non-null interior; by contrast, \cite{grun2008finite} provide a counter-example with a binary covariate. We satisfy this condition trivially --- for any covariate values --- because the misclassified covariate finite mixture in (\ref{mix.eq1}) defines mixture components with {\em parallel} hyperplanes, which share all parameters except the $\alpha_1 V$ term. Accordingly, component ordering is preserved across all covariate values. 

The identifiability of (\ref{mix.eq1}) ensures that parameter estimates can be obtained without side information on the magnitude of misclassification, such as validation data or known reclassification probabilities required by the existing methods. Instead, the specific form of the model being used allows estimation even though $V$ is not directly observed. Because this result may be counterintuitive, it is helpful to consider the normal response case with a dichotomous covariate and no concomitant variables. In that case, observed bimodality in the distribution of the response (for a fixed $V^{*}$) must be due to misclassification. Thus, simple visual inspection can be used to identify the misclassified observations, correct them, and proceed with the regression analysis. This simple process becomes much more difficult with concomitant variables, but a similar logic may be implemented through a formal likelihood-based analysis, which is the proposal of this paper.

\section{Asymptotic evaluation}\label{SecEst}
Based on the mixture representation in Equation (\ref{mix.eq1}), we can use the regular asymptotic theory for studying the efficiency loss due to the existence of covariate misclassification. 

\subsection{Fisher information}
Without loss of generality, we study the situation where there is no additional covariates $\mathbf{x}$ in the model. In this case, the likelihood for a single observation of $(Y,\,V^*)$ can be written as
\begin{equation}
L^*(\boldsymbol{\alpha},\,\boldsymbol{\varphi},\,\mathbf{Q},\,\boldsymbol{\pi}^*)\propto\left[\prod_{k=0}^{K-1}{(\pi_k^*)^{I_{(V*=k)}}}\right]\left[\sum_{j=0}^{K-1}q_{V^*j}\,f_Y(Y\,|\,\boldsymbol{\alpha},\,\boldsymbol{\varphi},\,V=j)\right],
\label{lik1}
\end{equation}
where $\boldsymbol{\pi}^*=\mathbf{P}'\boldsymbol{\pi}$ is the category probabilities for the observed variable $V^*$, $I_{(A)}$ is the indicator function on whether the event $A$ happens, and the reclassification probability only depends on $V^*$ and $j$.

Denote by $\boldsymbol{\theta}$ the vector that contains all the parameters in $\boldsymbol{\alpha}$, $\boldsymbol{\varphi}$, $\boldsymbol{\pi}^*$ and $\mathbf{Q}$. The corresponding expected Fisher information matrix has the form 
$\mathcal{I}^*(\boldsymbol{\theta})=E\left[\left.\left({\partial l^*(\boldsymbol{\theta})}/{\partial\boldsymbol{\theta}}\right)^{\otimes 2}\,\right|\,\boldsymbol{\theta}\right]$,
where $l^*(\boldsymbol{\theta})$ is the log-likelihood function for a single observation. 

For the asymptotic variance, we consider three \emph{successively harder} scenarios for estimating the parameters in $\boldsymbol{\theta}_{A}$, the block of parameters of interest including $\boldsymbol{\alpha}$ and $\boldsymbol{\varphi}$. 

The first case is \emph{when perfect observations of $(Y,\,V)$ are available} (i.e., when there is no misclassification). Denote by $l(\boldsymbol{\alpha},\,\boldsymbol{\varphi})=\log f_Y(Y\,|\,\boldsymbol{\alpha},\,\boldsymbol{\varphi},\,V)$ the log-likelihood function for a single observation of $(Y,\,V)$, and by $\mathcal{I}(\boldsymbol{\theta}_A)$ the corresponding expected Fisher information. Then the maximum likelihood estimator (MLE) $n\,\widehat{\boldsymbol{\theta}}_{A}\rightarrow MVN(n\,\boldsymbol{\theta}_{A},\,Acov_1)$, as $n\rightarrow\infty$. The asymptotic variance covariance matrix can be written as \begin{equation} Acov_0=[\mathcal{I}(\boldsymbol{\theta}_A)]^{-1}. \label{acov1} \end{equation}

The next easiest scenario is when the reclassification matrix $\mathbf{Q}$ is {\em known}, along with the observed values of $(Y,\,V^*)$. Because in practice the reclassification matrix is typically unknown, this represents an optimistic scenario; the methods of \cite{Shieh09} and \cite{buonaccorsi10} can be used in this case. For the MLE $\widehat{\boldsymbol{\theta}}_{A}$, the asymptotic variance covariance matrix is the corresponding block of covariance matrix in $[\mathcal{I}_{CC}^*(\boldsymbol{\theta})]^{-1}$. That is, 
\begin{equation} 
Acov_1=\left\{[\mathcal{I}^*_{CC}(\boldsymbol{\theta})]^{-1}\right\}_{AA},
\label{acov2}
\end{equation}
where the parameter block $\boldsymbol{\theta}_C$ consists of the parameters in $\boldsymbol{\theta}_A$ ($\boldsymbol{\alpha}$ and $\boldsymbol{\varphi}$) and $\boldsymbol{\theta}_{B}$ ($\boldsymbol{\pi}^*$), and $\left\{[\mathcal{I}^*_{CC}(\boldsymbol{\theta})]^{-1}\right\}_{AA} =
[\mathcal{I}^*_{AA}(\boldsymbol{\theta})]^{-1}- \mathcal{I}^*_{AB}(\boldsymbol{\theta}) [\mathcal{I}^*_{BB}(\boldsymbol{\theta})]^{-1} \mathcal{I}^*_{BA}(\boldsymbol{\theta})$.

Finally, the most challenging scenario is when we observe only $(Y,\,V^*)$, and the reclassification matrix $\mathbf{Q}$ is {\em unknown}. In this case, previous methods cannot be used but our new mixture approach can be used. When the mixture distribution is identifiable (see, e.g., \cite{yakowitz1968,atienza2007}), the MLE $\widehat{\boldsymbol{\theta}}_{A}$ will have consistency and asymptotic normality, with the asymptotic covariance matrix given by
\begin{equation} 
Acov_2=\left\{[\mathcal{I}^*(\boldsymbol{\theta})]^{-1}\right\}_{AA}.
\label{acov3}
\end{equation}

Given a parametric form for $f_Y(Y\,|\,\boldsymbol{\alpha},\,\boldsymbol{\varphi},\,V)$, the expected Fisher information can be evaluated numerically allowing to study the efficiency loss for inference on the parameters of interest, most importantly the regression coefficients corresponding to effects of interest. 

\subsection{Efficiency loss}
Here we illustrate the efficiency loss for a normal linear model when there is one binary covariate $V$ that is subject to nondifferential misclassification. Specifically, let $(Y\,|\,V)\sim N\left(\mu_{V}, \sigma^2\right)$, with the conditional mean $\mu_V=\alpha_0+\alpha_1V$, and constant conditional variance, $\sigma^2$.


Denote by $\widehat{\alpha}_1$ the MLE of the covariate effect $\alpha_1$, and $\boldsymbol{\theta}=(\alpha_0,\,\alpha_1,\,\sigma^2,\,\boldsymbol{\pi}^*,\,\mathbf{Q})$. For the easiest scenario, when $(Y,\,V)$ are observed, standard asymptotic theory gives $n\,\mbox{Var}(\widehat{\alpha}_1)\rightarrow Avar_0\left(\boldsymbol{\theta}\right)$, with a closed-form $Avar_0\left(\boldsymbol{\theta}\right)={\sigma^2}/{[\pi_1(1-\pi_1)]}$. In order to obtain the forms of the asymptotic variance terms in Equations (\ref{acov2}) and (\ref{acov3}), we can write the log-likelihood for a single sample of $(Y,\,V^*)$ as
\begin{align}
l^*(\boldsymbol{\theta})=&\log\left[q_{V^*0}\,f(Y\,|\,\alpha_0,\sigma^2)+(1-q_{V^*0})\,f(Y\,|\,\alpha_0+\alpha_1,\sigma^2)\right]\notag\\
&+{V^*}\log(\pi^*_1)+({1-V^*})\log(1-\pi^*_1)+C,
\end{align}
where $C$ is a constant, and $f(Y\,|\,\alpha_0,\sigma^2)$ denotes the normal density function with mean $\alpha_0$ and variance $\sigma^2$. 

Based on the log-likelihood function, we can obtain the expected Fisher information matrix as
\[\mathcal{I}^{*}(\boldsymbol{\theta})=\pi^*_1E[s^{\otimes 2}(\boldsymbol{\theta})|V^*=1]+(1-\pi^*_1)E[s^{\otimes 2}(\boldsymbol{\theta})|V^*=0],\]
where the score vector is given by
\begin{equation}
s(\boldsymbol{\theta})=\left(\frac{\partial\,l^*(\boldsymbol{\theta})}{\partial\,\alpha_0},\,\frac{\partial\,l^*(\boldsymbol{\theta})}{\partial\,\alpha_1},\,\frac{\partial\,l^*(\boldsymbol{\theta})}{\partial\,\sigma},\,\frac{\partial\,l^*(\boldsymbol{\theta})}{\partial\,\pi^*_1},\,\frac{\partial\,l^*(\boldsymbol{\theta})}{\partial\,q_{00}},\,\frac{l^*(\boldsymbol{\theta})}{\partial\,q_{10}}\right)',
\end{equation}
with the analytical forms of the elements given in Section A in the supplement (\cite{XiaHahn16sup}).

Using numerical integration methods to evaluate each of the elements in $E[s^{\otimes 2}(\boldsymbol{\theta})|V^*=v^*]$, it is possible to compute the asymptotic variance terms in Equations (\ref{acov2}) and (\ref{acov3}). Numerical integration is more stable after converting the infinite integral to a finite integral by transforming the variables using their cumulative distribution functions. We use the R function $integral()$ from the \textbf{pracma} package, and found the adaptive methods based on the Gauss–Kronrod quadrature (\cite{davis2007}) to be fast and reliable. 

\subsection{Numerical evaluation}
Here, we study the {\em ratio of the asymptotic standard deviations}, denoted $Rasd_j=\sqrt{Avar_j/Avar_0}$ ($j=1,\,2$), respectively for the two cases when the reclassification probabilities $q_{01}$ and $q_{10}$ are known and unknown. For the normal model, the ratio of the asymptotic standard deviations $Rasd_j$ only depends on three factors: the effect size $\alpha_1/\sigma$, the severity of misclassification, and the binomial proportion $\pi_1$. In Figure \ref{eff1} and Web Figures 1 through 2 in the supplement, we investigate five cases: $\alpha_1=0.3$, 0.5, 1, 2 and 5 with fixed $\sigma=1$. In Figure \ref{eff1}, the binomial proportion $\pi_1 = 0.5$ (Web Figure 2 in the supplement presents the results for $\pi_1 = 0.2$, \cite{XiaHahn16sup}).

\begin{figure}[htp]
   \centering
      \mbox{\subfigure[$Rasd_1:\,\alpha_1=1,\,\pi_1=0.5$]{\epsfig{figure=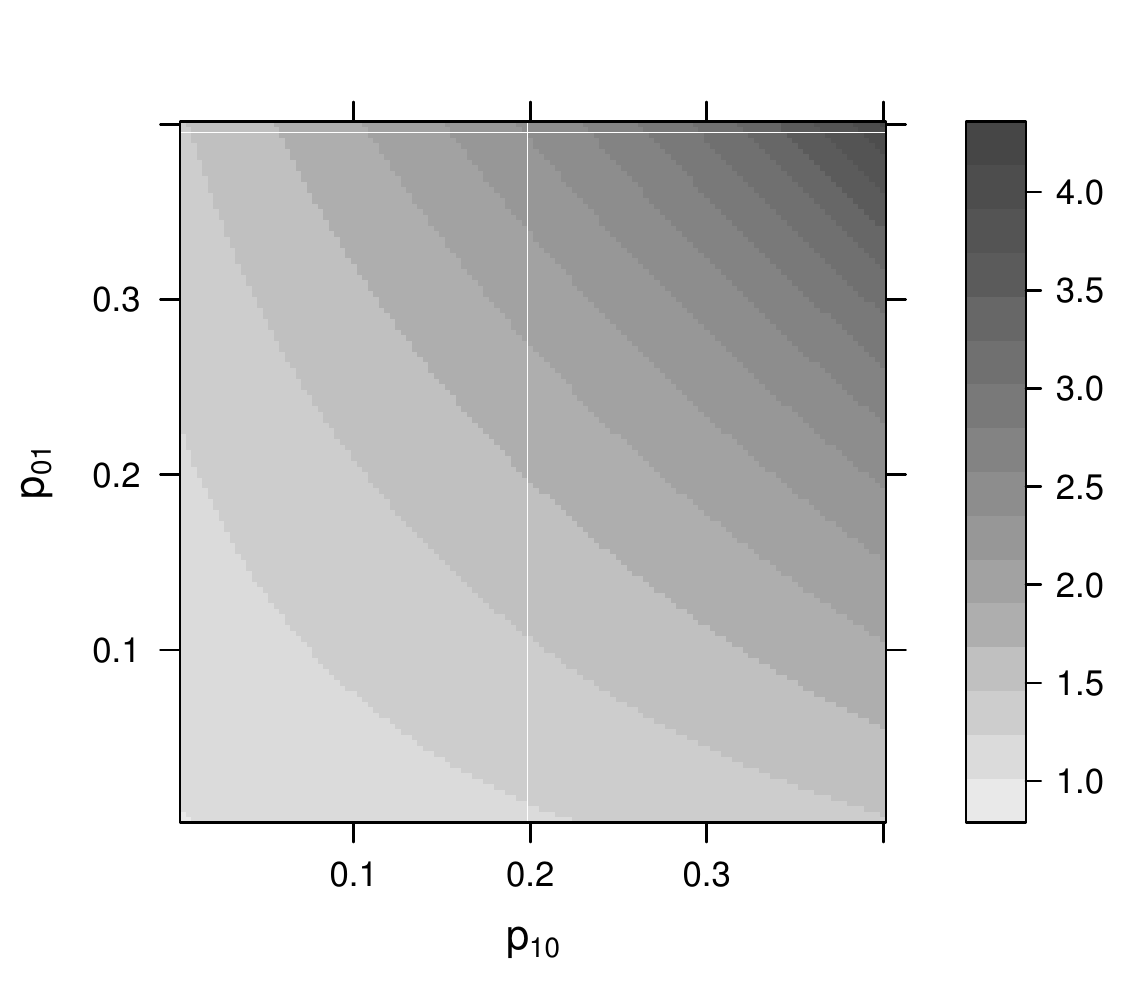,height=1.9in}}\quad
         \subfigure[$Rasd_1:\,\alpha_1=2,\,\pi_1=0.5$]{\epsfig{figure=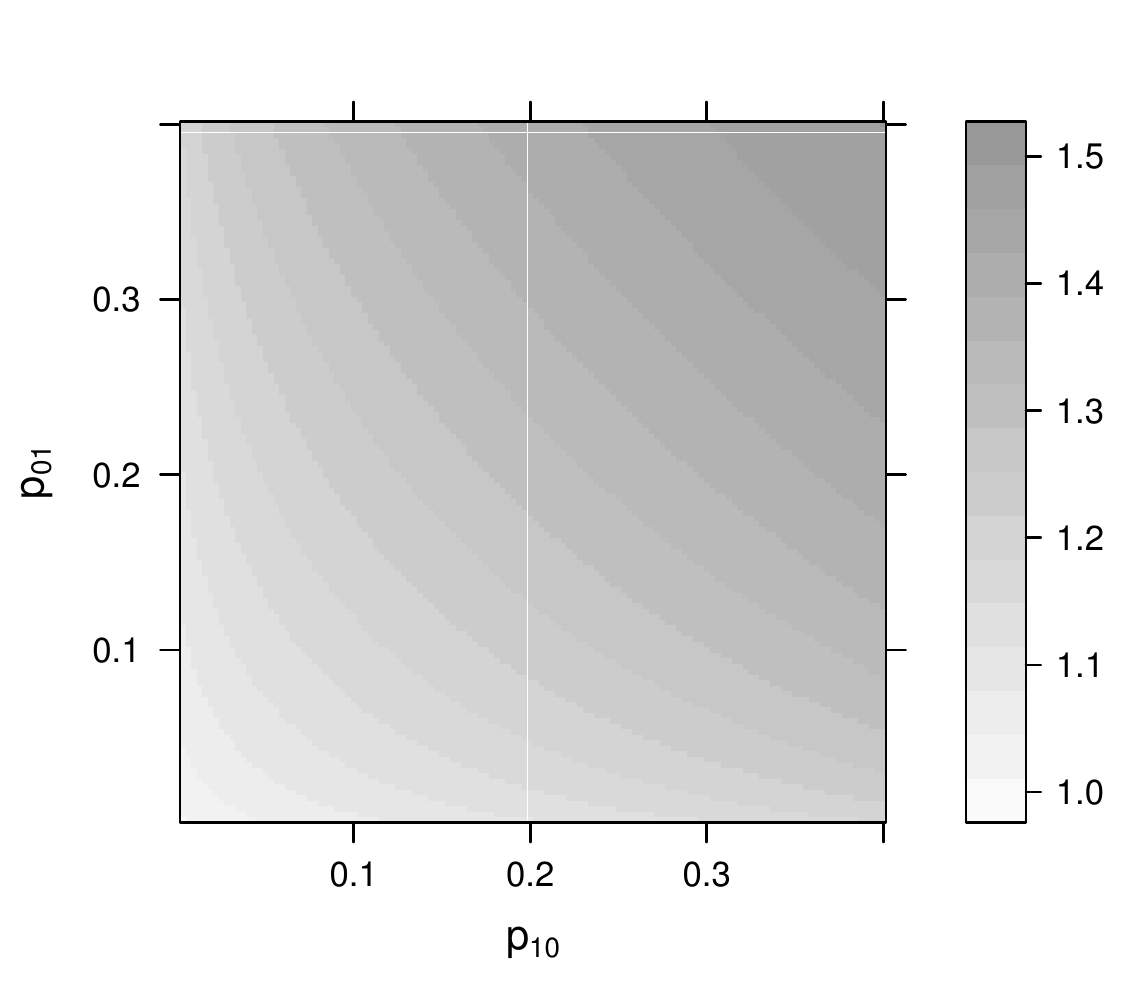,height=1.9in}}\quad
         \subfigure[$Rasd_1:\,\alpha_1=5,\,\pi_1=0.5$]{\epsfig{figure=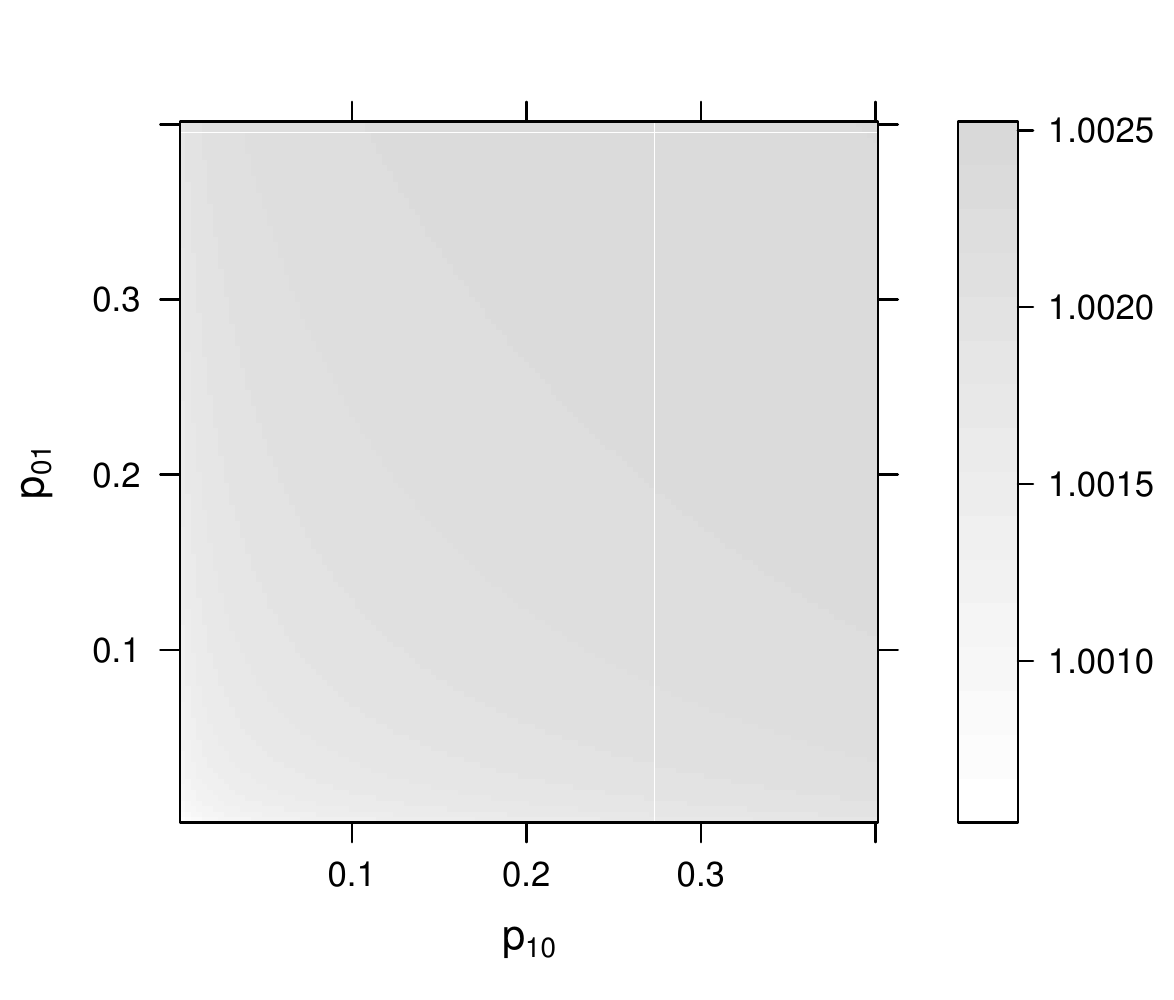,height=1.9in}}}  \\
     \mbox{\subfigure[$Rasd_2:\,\alpha_1=1,\,\pi_1=0.5$]{\epsfig{figure=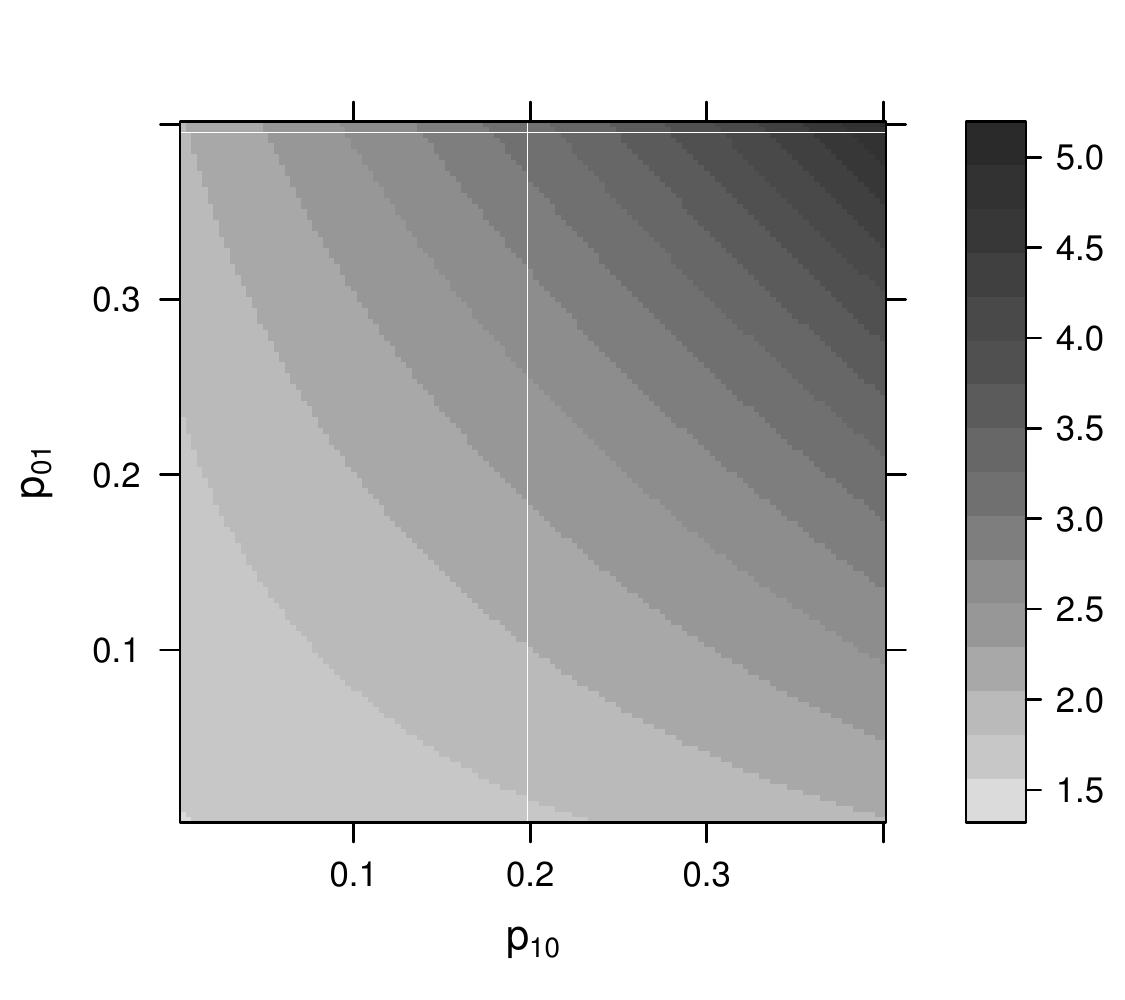,height=1.9in}}\quad
         \subfigure[$Rasd_2:\,\alpha_1=2,\,\pi_1=0.5$]{\epsfig{figure=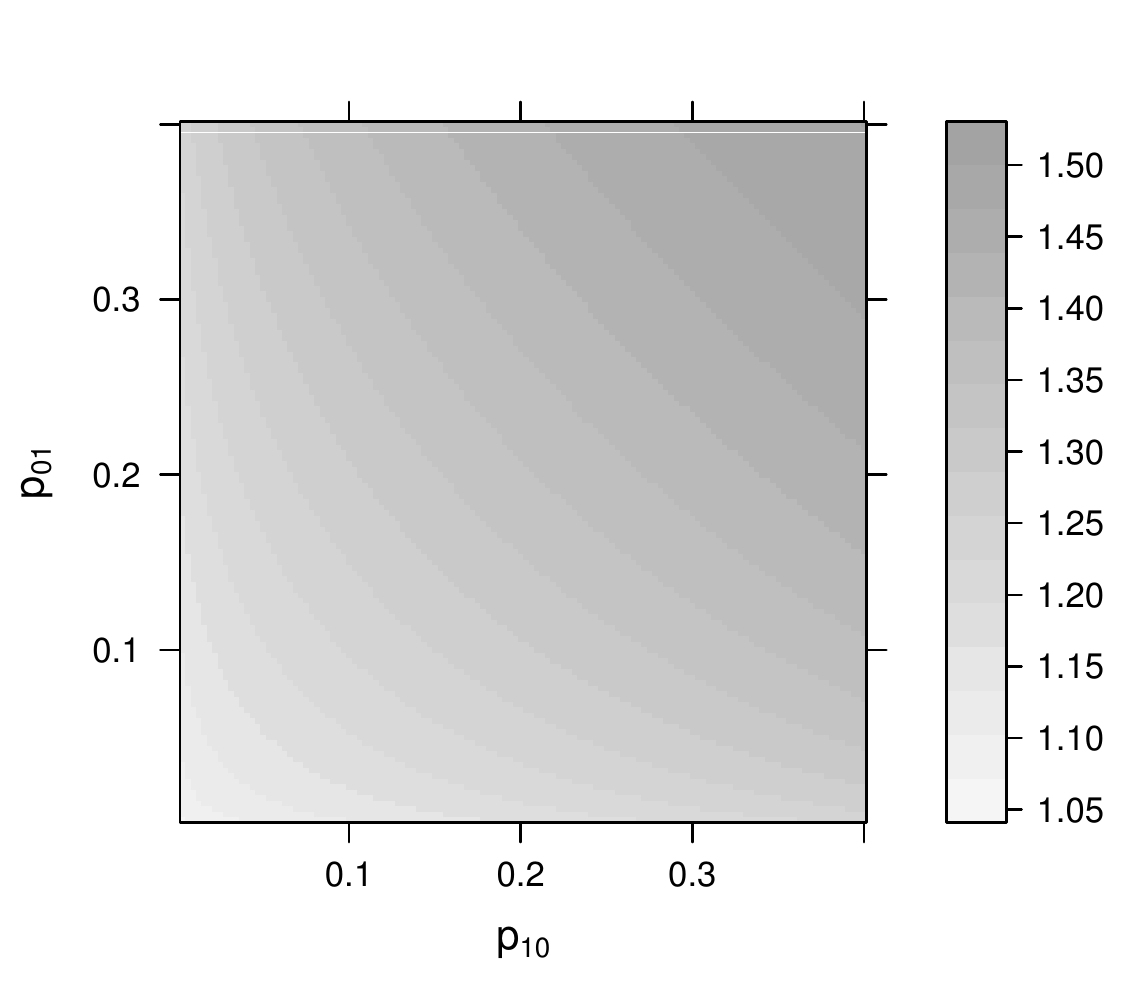,height=1.9in}}\quad
         \subfigure[$Rasd_2:\,\alpha_1=5,\,\pi_1=0.5$]{\epsfig{figure=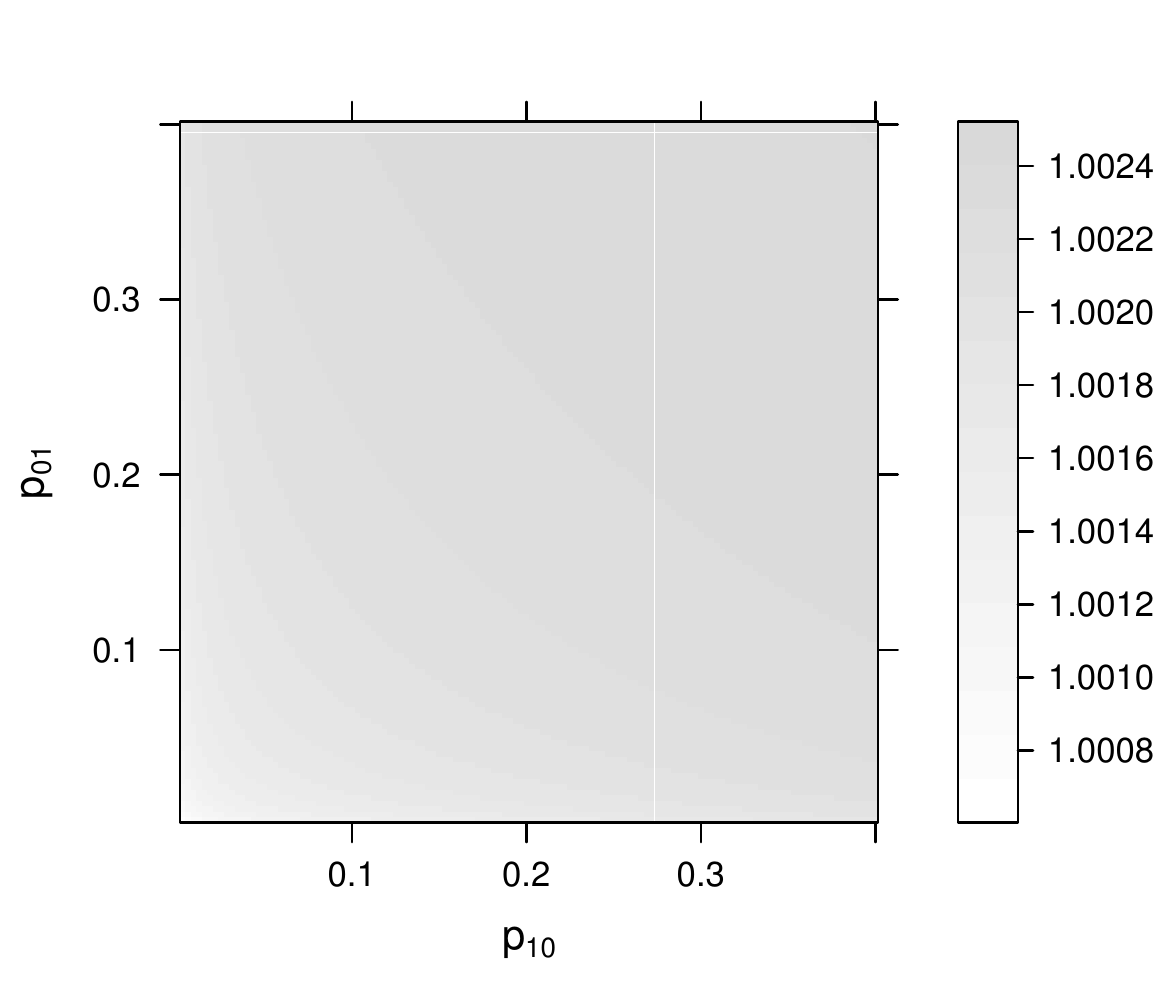,height=1.9in}}}    
     \caption{Efficiency loss for the effect $\alpha_1$ with varying misclassification probabilities, under the normal model. The top panels present $Rasd_1$, the ratio of the asymptotic standard deviations for the case when the reclassification matrix is known versus that where the true status $V$ is being observed. The bottom panels correspond to $Rasd_2$, that between the proposed model without knowing the reclassification matrix and the true model with accurately measured $V$.}
   \label{eff1}
\end{figure}

Figure \ref{eff1} shows that the efficiency loss is dominated by the effect size and the severity of the misclassification. For most of the scenarios, the ranges of the efficiency loss across different combinations of misclassification probabilities are similar, for the two cases when the reclassification matrix is known and unknown. Knowing the reclassification matrix offers more benefit in the cases when the effect size is small. When the effect size is large (i.e., $\geq 2$), there is very little efficiency loss, regardless of whether the reclassification probabilities are known or not. A smaller or larger binomial proportion $\pi_1$ will result in mildly larger benefit from knowing the reclassification probabilities. In summary, large efficiency loss only occurs when the misclassification probabilities are large and the effect size is small. In these cases, the quality of the measurement is so questionable that a much larger confidence interval is needed, even for cases when we know the reclassification matrix. Note that on the boundaries of the figures, when one of the misclassification probabilities is 0 or 1, the ratio of the asymptotic standard deviations $Rasd_2$ will be infinite, owing to the nonidentifiability of mixture models in these scenarios (see, e.g., \cite{li2007, Xia16}).

\section{Simulation studies} \label{Simu}
In order to illustrate the performance of the method in finite samples, we perform Monte Carlo simulation studies. Inference is performed using a fully Bayesian approach, based on Markov chain Monte Carlo (MCMC) samples from the posterior. This approach provides convenient statistical inference and the capability of incorporating prior information on the reclassification probabilities. For frequentist estimation, the {\tt R} package \textbf{FlexMix} (\cite{grun2008}) can be used to mixture models using the EM algorithm, including the specific mixture models given in Equations (\ref{mix.eq1}) and (\ref{mix.eq2}).

Here, to illustrate our method, we consider two cases: the normal response case and the Poisson response case. More comprehensive studies for the normal, Poisson and also gamma models are given in the supplement (\cite{XiaHahn16sup}). 

\subsection{Normal model}
\subsubsection{Impact from effect size and external information}
Here, we investigate the impact of the effect size and the (in)accuracy of external information concerning the reclassification probabilities on posterior estimates. In the normal linear model, we illustrate the convergence of the credible interval when there is an ordinal covariate $V$ that is subject to nondifferential misclassification. In particular, we specify the conditional distribution of the response variable as $(Y\,|\,V)\sim N\left(\mu_{V}, \sigma^2\right)$, with the conditional mean having a linear form $\mu_V=\alpha_0+\alpha_1V$. 

We assume that the ordinal covariate has a multinomial distribution, with probabilities $\boldsymbol{\pi}=(\pi_0,\pi_1,\pi_2)=(0.2,0.3,0.5)$ corresponding to the three values 0, 1 and 2. The following classification and corresponding reclassification matrices are assumed for obtaining the misclassified sample of $V^*$. 
\[\mathbf{P}=
\left[
\begin{array}{lll}
0.80&  0.15&  0.05\\
0.10&  0.70&  0.20\\
0.05&  0.15&  0.80\\
\end{array}
\right]\quad\mbox{or}\quad
\mathbf{Q}=
\left[
\begin{array}{lll}
 0.74& 0.14& 0.12\\
 0.10& 0.67& 0.24\\
 0.02& 0.13& 0.85\\
\end{array}
\right].\]
After generating the sample of $V$ from its multinomial distribution, we generate $Y$ with variance $\sigma^2=4$, and regression coefficients $(\alpha_0,\,\alpha_1)$ being either $(12, 2)$, $(12, 4)$ or (12, 10). These three sets of parameters correspond to the same effect sizes as those in Figure \ref{eff1} (i.e., $\alpha_1 / \sigma$, of 1, 2 and 5, respectively).

We illustrate the finite sample performance of four models with the sample sizes of 100, 400, 1,600, 6,400 and 25,600. The {\em naive model} refers to linear regression using the observed values of $V^*$ as if there were no misclassification. The {\em true model} refers to linear regression using the correct classification $V$. The third model assumes that the reclassification matrix $\mathbf{Q}$ is known. Finally,  our finite mixture model estimates the reclassification matrix based on the mixture representation. 

For all four models, a normal prior with mean 0 and variance 100 is used for the intercept, and a gamma prior with parameters 0.0001 and 0.1 is used for the slope. The gamma prior represents our knowledge on the sign of the effect obtained from the naive estimate (i.e., the assumption that misclassification will not change the direction of the effect, in order to ensure the identifiability of the mixture model). For the normal model with an ordinal covariate, we place vague Dirichlet priors on the probability parameters in $\boldsymbol{\pi}^*$and each column of $\mathbf{Q}$, with concentration parameters of 1. The prior mean and standard deviation of the Dirichlet prior are 1/3 and 0.24. The above Dirichlet prior represents a situation when we have no external information on these probabilities. For the reclassification probabilities in the cases with $\alpha_1/\sigma=1$ and 5, we consider an additional set of informative priors for the proposed model. We use $60\times q_{kj}$ as the concentration parameter for each of the reclassification probabilities, in order to match the prior mean with the true value. The resulting priors have small standard deviation ranging from 0.02 to 0.06. For the case where $\alpha_1/\sigma=2$, we study the model performance when the external information the reclassification probabilities is inaccurate. Instead of knowing the true values, we assume the magnitudes of misclassification are under-estimated. In particular, we assume the estimates from external data are
\[\mathbf{\widehat{Q}}=
\left[
\begin{array}{lll}
 0.84& 0.04& 0.12\\
 0.10& 0.77& 0.14\\
 0.02& 0.03& 0.95\\
\end{array}
\right].\]
Thus, we use the above estimates in the third model assuming known reclassification probabilities. For the proposed model, we consider an additional set of Dirichlet priors with the corresponding concentration parameters given by $60\times\mathbf{\widehat{Q}}$ (i.e., representing an inaccurate guess). We take every 10th sample for all the models, after dropping the first 15,000 samples as burn-in. For a posterior sample of 5,000, the effective sample size is over 4,500 for all the models. The 95\% equal-tailed credible intervals of the regression effect and other parameters from the four models are provided in Figure \ref{normal1}, and Web Figures 3 through 10 in the supplement. 

\begin{figure}[htp]
   \centering
   \mbox{\subfigure[$\alpha_1=2$ with accurate $\mathbf{\widehat{Q}}$]{\epsfig{figure=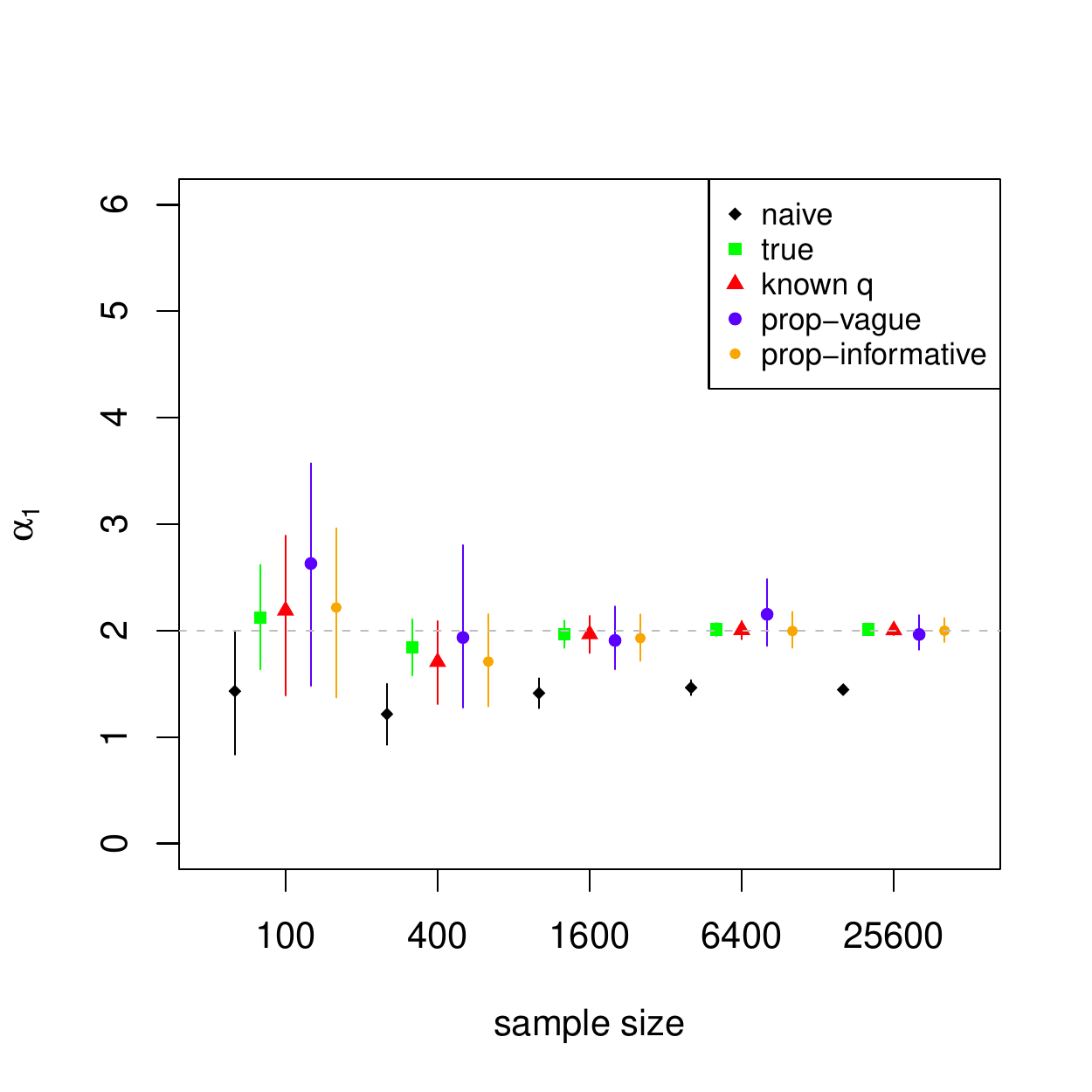,height=2.2in}}\quad
	       \subfigure[$\alpha_1=4$ with biased $\mathbf{\widehat{Q}}$]{\epsfig{figure=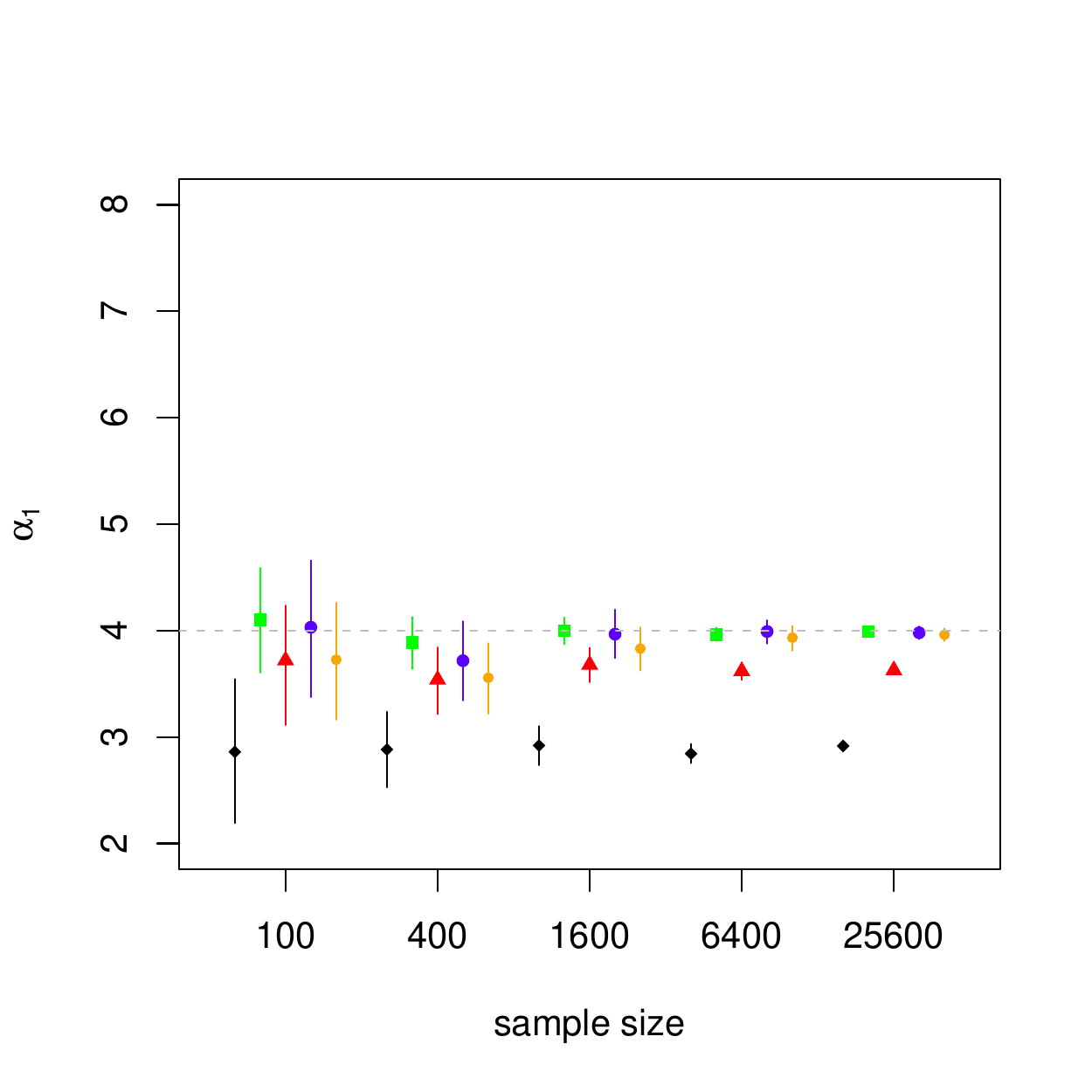,height=2.2in}}\quad
         \subfigure[$\alpha_1=10$ with accurate $\mathbf{\widehat{Q}}$]{\epsfig{figure=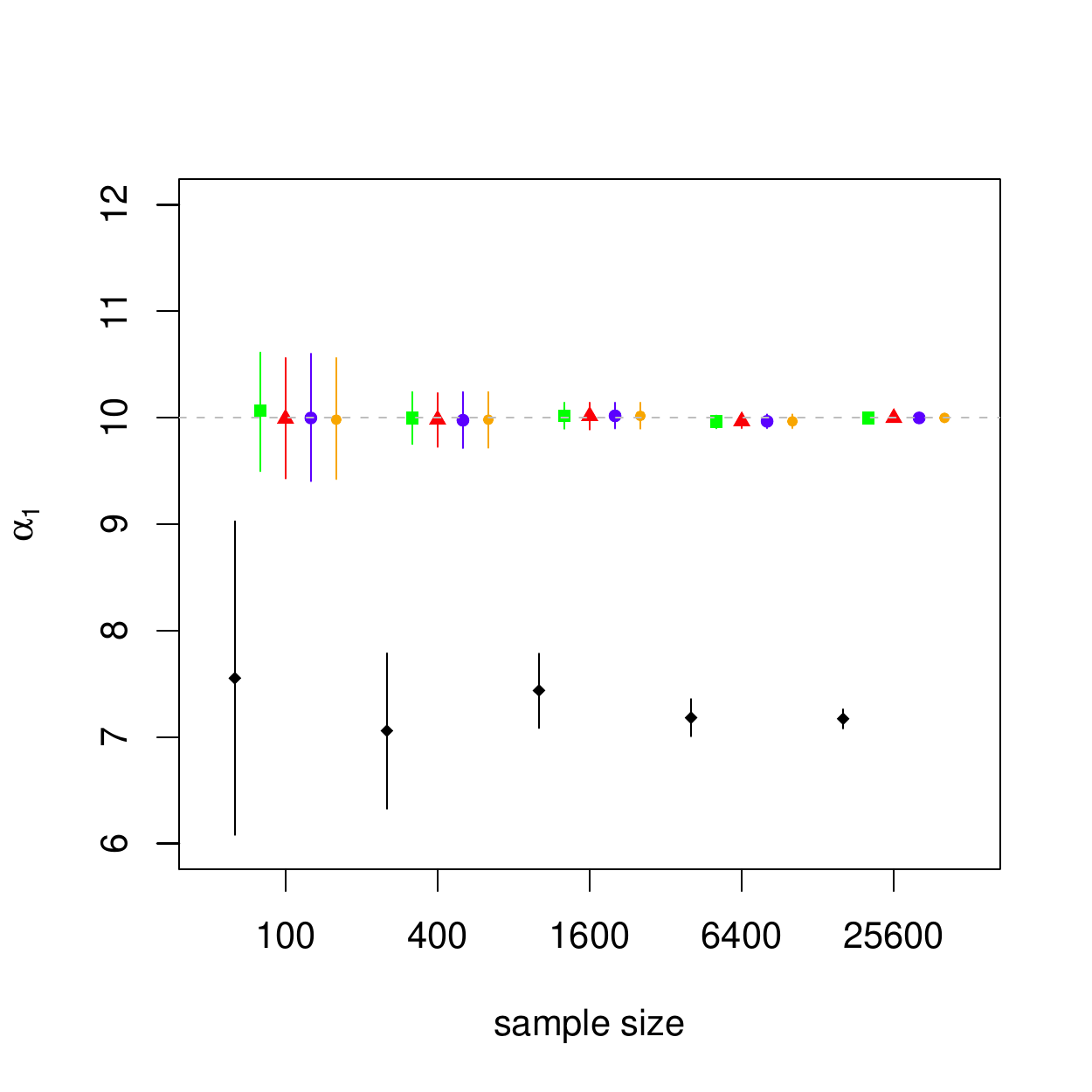,height=2.2in}}}     
     \caption{Equal-tailed credible intervals of $\alpha_1$ for normal responses with an ordinal $V$. With $\sigma^2=4$, the values of $\alpha_1$ in the three panels correspond to an effect size of 1, 2 and 5. In the first panel, observe that the informative prior (gold) gives a posterior interval almost identical to the fixed $\mathbf{\widehat{Q}}$ case (red), both of which have notably smaller variance than the vague prior case, especially at lower sample sizes. In the second panel, observe that the informative prior is able to ``unlearn'' its prior bias when $\mathbf{\widehat{Q}}$ is incorrect, thus outperforming the fixed $\mathbf{\widehat{Q}}$ method (red) for sample sizes larger than 400. The final panel shows that when the effect size is large enough, all methods that account for misclassification are comparable with the case where the true status is observed.}
   \label{normal1}
\end{figure}

Figure \ref{normal1} reveals several intuitive patterns. First, the third panel shows that when the effect size is large, observing the true covariate status or knowing the reclassification probabilities confers no benefit for parameter inference --- the Bayesian mixture analysis under a vague prior gives equivalent results. Second, the first panel shows, that informative priors give performance very similar to method which assumes known reclassification probabilities in terms of point estimate and posterior variance, and, by contrast, the vague prior has much wider posterior intervals; this shows that prior information can be incorporated successfully via a prior.  Finally, the second panel shows that when the reclassification probabilities are incorrectly specified, an informative prior (centered at the incorrect value) can eventually overcome this bias, outperforming the method which fixes the reclassification probabilities at the wrong value; incorporating information about the reclassification probabilities via a prior is therefore a more conservative approach.  

As expected, there is large bias in the naive estimates (which ignore misclassification). For the normal model, the variance can be dramatically over-estimated when misclassification is ignored. In the supplement (\cite{XiaHahn16sup}), it is shown that the reclassification probabilities can also be learned from the observed response and covariate, without further information. The learning is more efficient in the case when the effect size is large.

\subsubsection{Impact from tail heaviness}
For the normal model, we perform a sensitivity analysis on the impact from model misspecification. In particular, we consider scenarios where the response variables are generated from a Student-$t$ distribution with the location and scale parameters denoted by $\mu_V$ and $\sigma$. We assume the same normal model as the previous subsection, and consider the large-effect case with $\alpha_0=12$ and $\alpha_1=10$. The category and reclassification probabilities also take the same values. 

We consider three cases in an attempt to understand the impact of tail heaviness, as well as that of larger variability caused by either a larger scale parameter or a smaller number of degrees of freedom for the Student-$t$ distribution. For the first scenario, we assume that the number of degrees of freedom $\nu=20$. We choose the scale parameter to be $\sigma=\sqrt{3.6}$, so that the variance is the same (4) as the large-effect case studied in the previous subsection. For the second and third scenarios, we increase the variance of the Student-$t$ distribution by increasing the scale parameter $\sigma$ to $\sqrt{6.48}$ or decreasing the number of degrees of freedom $\nu$ to 4. Both cases correspond to a variance of 7.2.

We consider five sample sizes, as in the previous subsection. We keep the MCMC and prior settings the same as the earlier case. The 95\% equal-tailed credible intervals of the regression effect and other parameters are provided in Figures \ref{normal.sen1}, and Web Figures 11 through 16 in the supplement. 

\begin{figure}[htp]
   \centering
   \mbox{\subfigure[$\sigma^2=3.6$, $\nu=20$]{\epsfig{figure=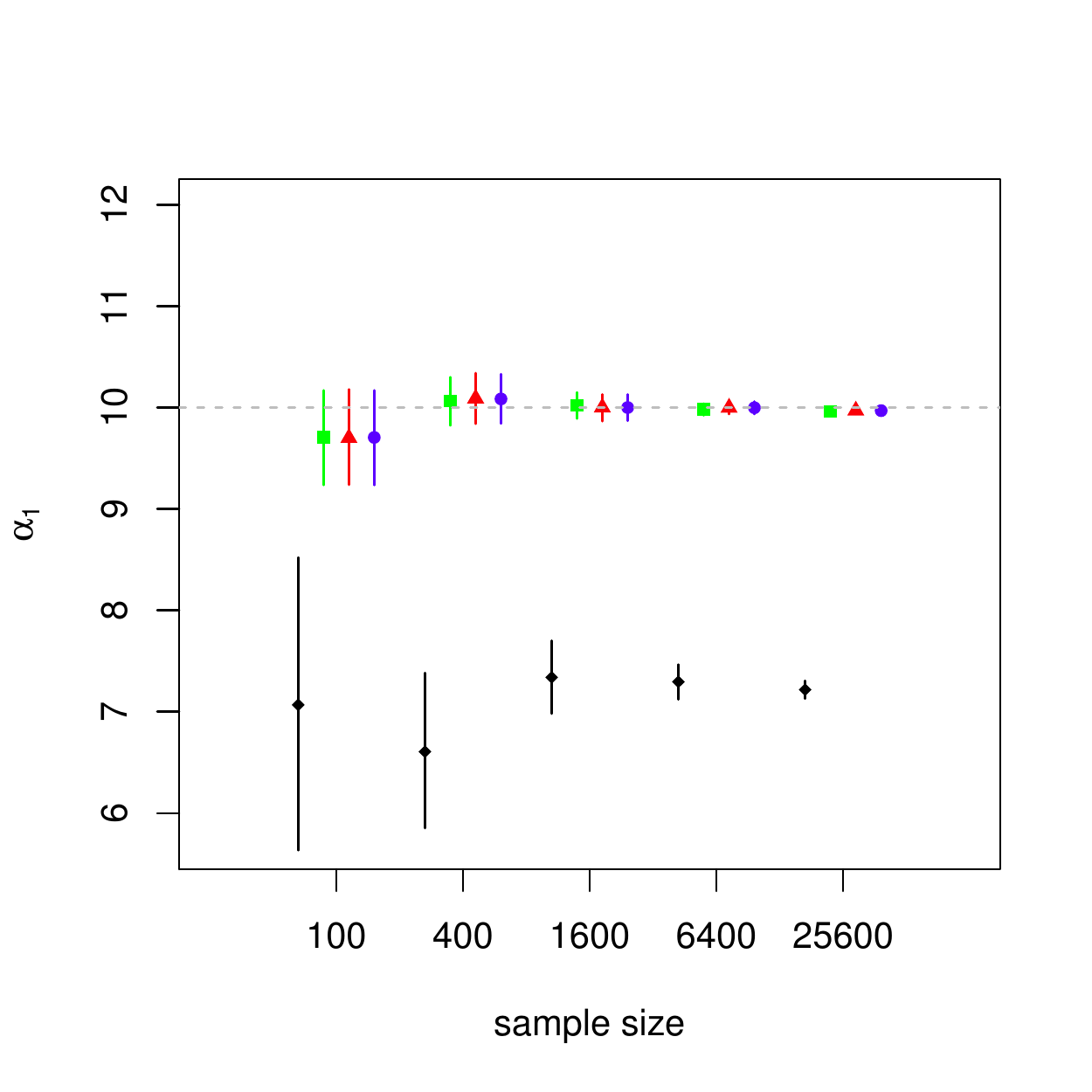,height=2.2in}}\quad
         \subfigure[$\sigma^2=6.48$, $\nu=20$]{\epsfig{figure=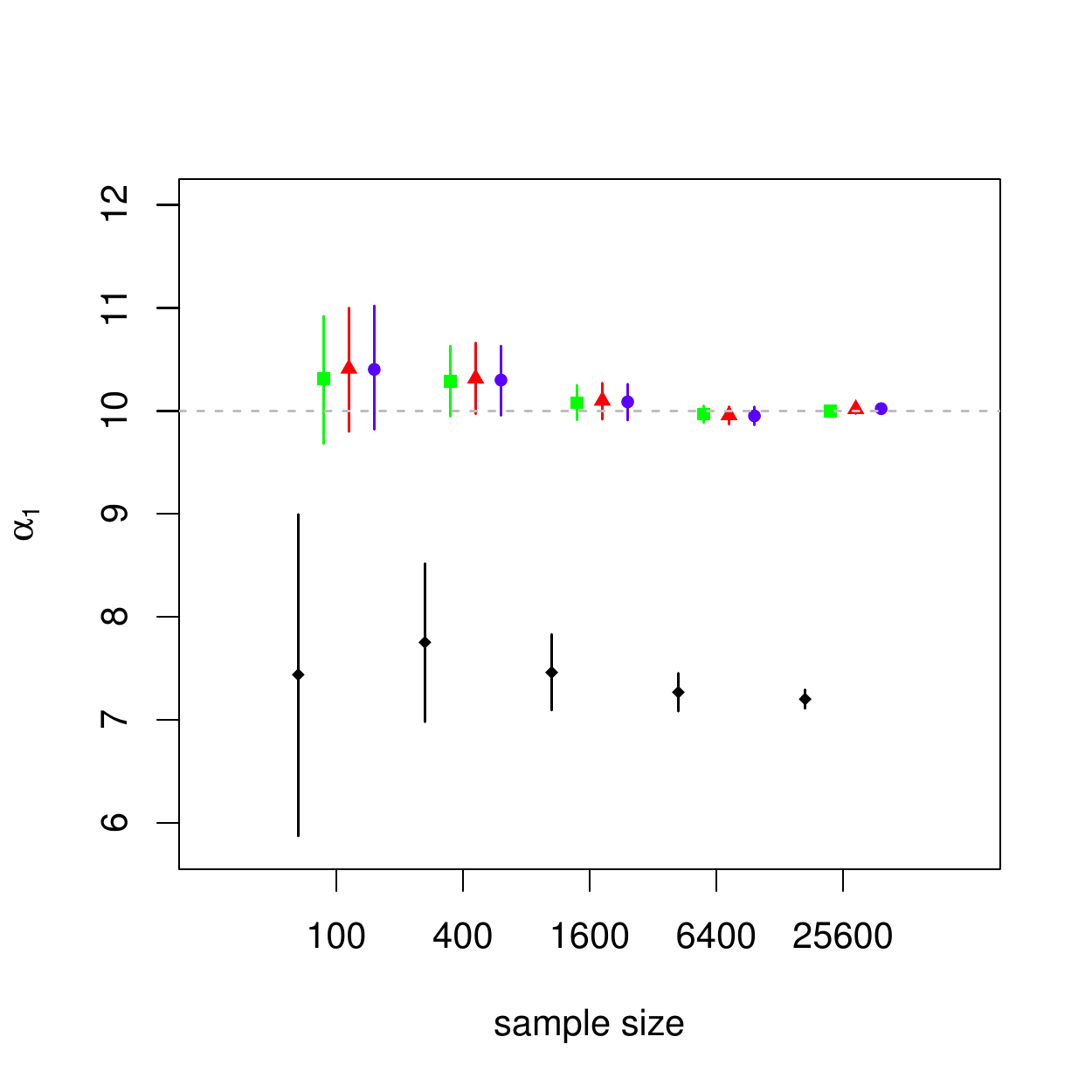,height=2.2in}}\quad
         \subfigure[$\sigma^2=3.6$, $\nu=4$]{\epsfig{figure=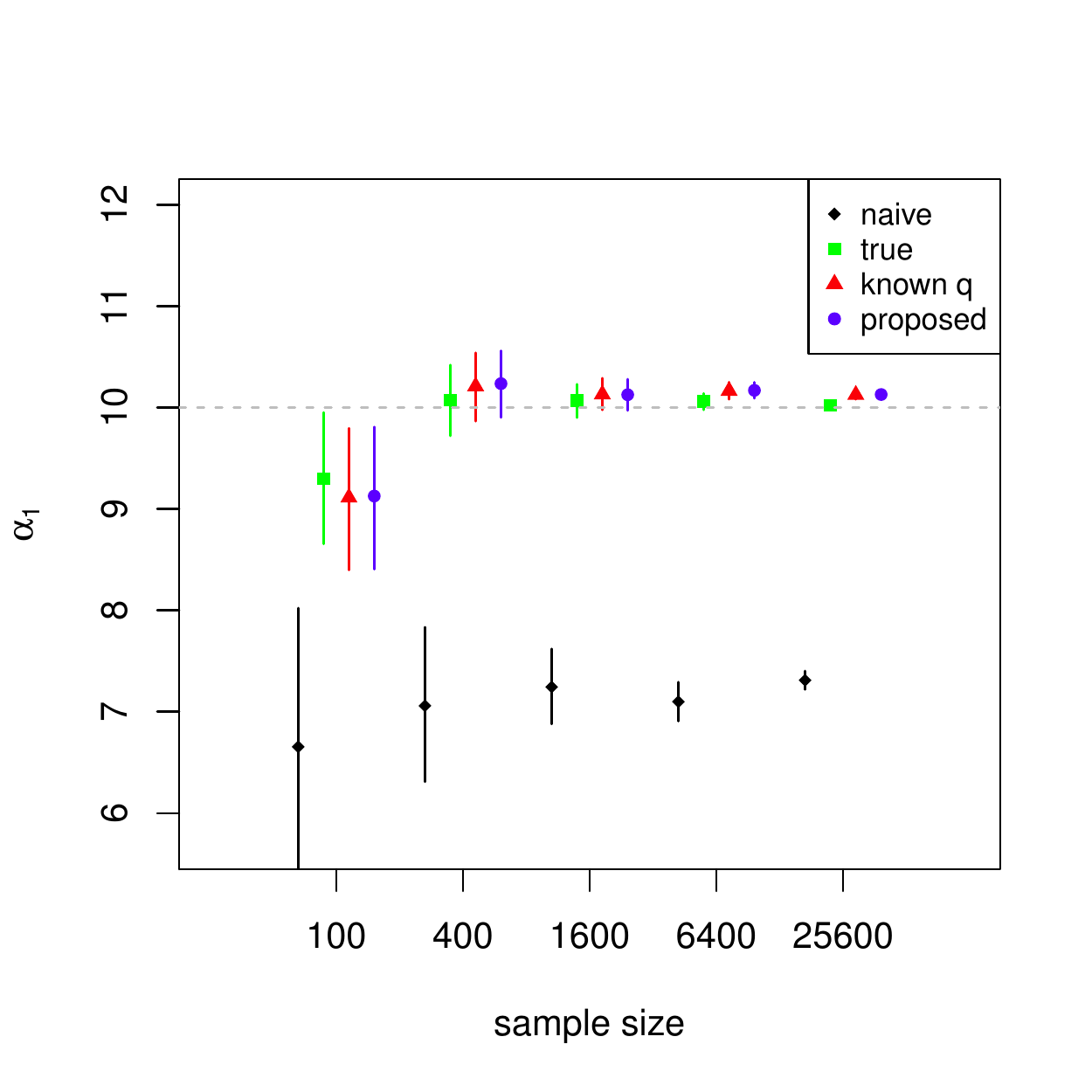,height=2.2in}}}     
     \caption{Equal-tailed credible intervals of $\alpha_1$ for Student t responses with an ordinal $V$. The variances for the three cases are 4, 7.2 and 7.2. With a slope $\alpha_1=10$, these sets of values correspond to an effect size of 5, 3.7 and 3.7, respectively. }
   \label{normal.sen1}
\end{figure}

We observe that: (1) mild tail heaviness does not seem to have an impact on the performance of the normal model; (2) increasing the scale parameter alone seems to lead to larger variability in the estimation (i.e., owing to a smaller effect size); (3) increasing the tail heaviness results in both small biases in the estimation and larger variability. The impact seems to be the same, regardless of whether the reclassification probabilities are known or not. The bias in the estimation is probably caused by the fact that lower (upper) tail values from both mixture components are more likely to be attributed to the group with a lower (higher) mean, when the tails are very heavy. 

\subsection{Poisson model}
\subsubsection{Impact from severity of misclassification}
For the Poisson model, we study the impact of the severity of misclassification on finite sample inference. We assume there is a binary covariate $V$ that is subject to nondifferential misclassification. We generate the sample of $V$ using a Bernoulli trial with the probability $\pi_1=0.5$. In particular, we specify the conditional distribution of the response variable as $(Y\,|\,V)\sim Poisson\left(\mu_{V}\right)$, with the conditional mean $\mu_V=\exp(\alpha_0+\alpha_1V)$. The corresponding sample of Poisson counts $Y$ is then generated from that of $V$, with the regression coefficients $(\alpha_0,\,\alpha_1)$ assumed to be $(1.2, 1)$. Two different scenarios, $(p_{01},p_{10})=(0.1,0.125)$ and $(p_{01},p_{10})=(0.125,0.25)$, are assumed for the misclassification probabilities for obtaining the corresponding sample of $V^*$. 

Similar to the case of the normal model, we explore the finite sample performance at the sample sizes of 100, 400, 1,600, 6,400 and 25,600. We compare the equal-tailed credible intervals from the naive model, the true model, the model with the known reclassification matrix $\mathbf{Q}$, and the proposed model as defined in the earlier subsection. For all models, independent normal priors with mean 0 and variance 10 are used for the regression coefficients $\alpha_0$ and $\alpha_1$. For the parameters $p_{01}$, $p_{10}$ and $\pi_0$, independent uniform priors on $(0, 1)$ are used. The MCMC details are the same as the normal case. The 95\% equal-tailed credible intervals of the regression effect and other parameters from the four models are provided in Figure \ref{poisson1}, and Web Figures 17 and 19 in the supplement. 
\begin{figure}[htp]
   \centering
   \mbox{\subfigure[$p_{10}=0.125$]{\epsfig{figure=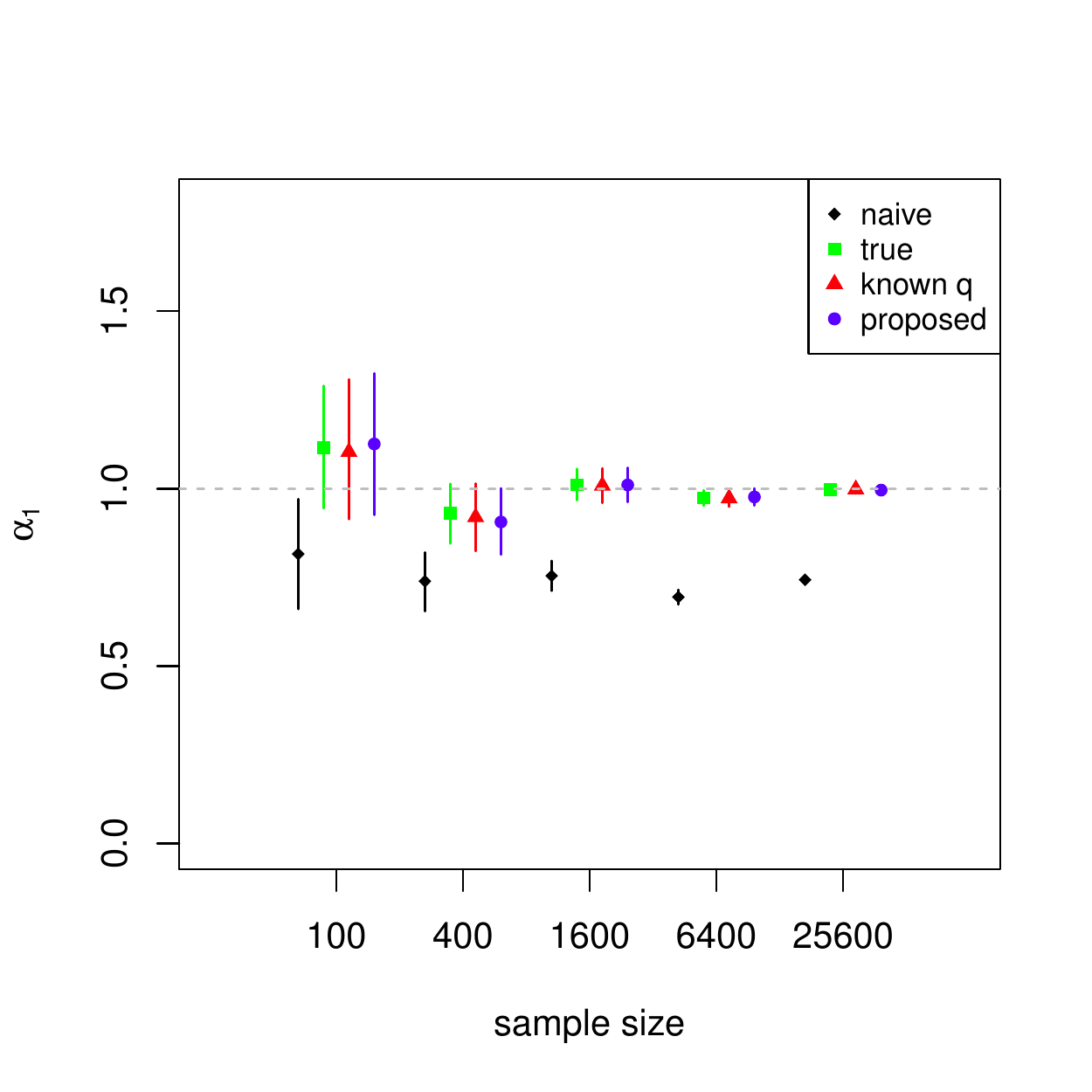,height=2.2in}}\quad
         \subfigure[$p_{10}=0.25$]{\epsfig{figure=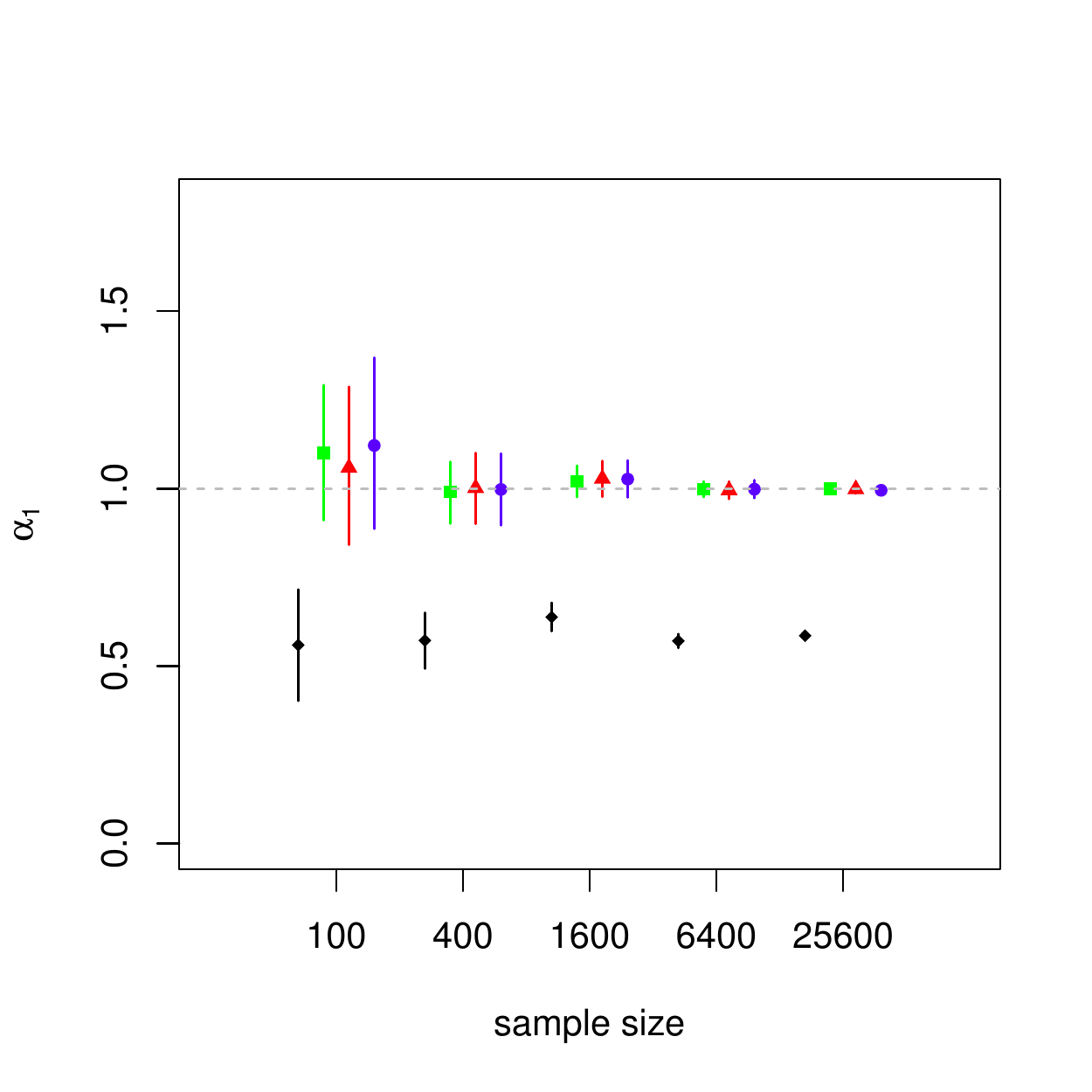,height=2.2in}}}     
     \caption{Equal-tailed credible intervals of $\alpha_1$ for Poisson responses with a binary $V$. }
   \label{poisson1}
\end{figure}

The figures show that the naive model demonstrates larger attenuation bias in the second scenario, when the misclassification probabilities are larger. Because the effect size is large, there is little efficiency loss regardless of the severity of the misclassification and knowing the reclassification matrix offers little extra efficiency. For both cases, the credible intervals of the misclassification probabilities, $p_{01}$ and $p_{10}$, are seen to converge to the true values with increasing sample size. Knowing the reclassification probabilities helps with the efficiency of estimating their misclassification counterparts. These results confirm the conclusions from an asymptotic evaluation based on the observed Fisher information. 

\subsubsection{Impact from zero inflation}
For count data, it will be interesting to perform a sensitivity analysis on the impact from model misspecification when there is zero inflation. Here, we assume that the response variable $Y$ is generated from a zero-inflated Poisson (ZIP) distribution. We assume that the conditional mean of the ZIP distribution is given by $\mu'_V=(1-w)\mu_V=\exp(\alpha_0+\alpha_1V)$, with $w$ being the percentage of additional zeros and $\mu_V$ being the mean of the Poisson component. For the sensitivity analysis, we take the less-severe misclassification scenario with $(p_{01},p_{10})=(0.1,0.125)$. We study two scenarios on zero inflation with $w=5\%$ and $10\%$, a fixed percentage of zeros across the $V=0$ and $V=1$ groups. We consider two cases with $(\alpha_0,\,\alpha_1)=(1.2,\,1)$ and (0, 1), with the conditional distributions of $Y$ presented in Web Figure 20. 

Again, we study the performance of the four models under five sample sizes. The prior and MCMC settings are the same as those in the previous subsection. The 95\% equal-tailed credible intervals of the regression effect and other parameters are provided in Figures \ref{poisson.sen1}, and Web Figures 21 through 24 in the supplement. 

\begin{figure}[htp]
   \centering
   \mbox{\subfigure[$\alpha_0=1.2$, $w=5\%$]{\epsfig{figure=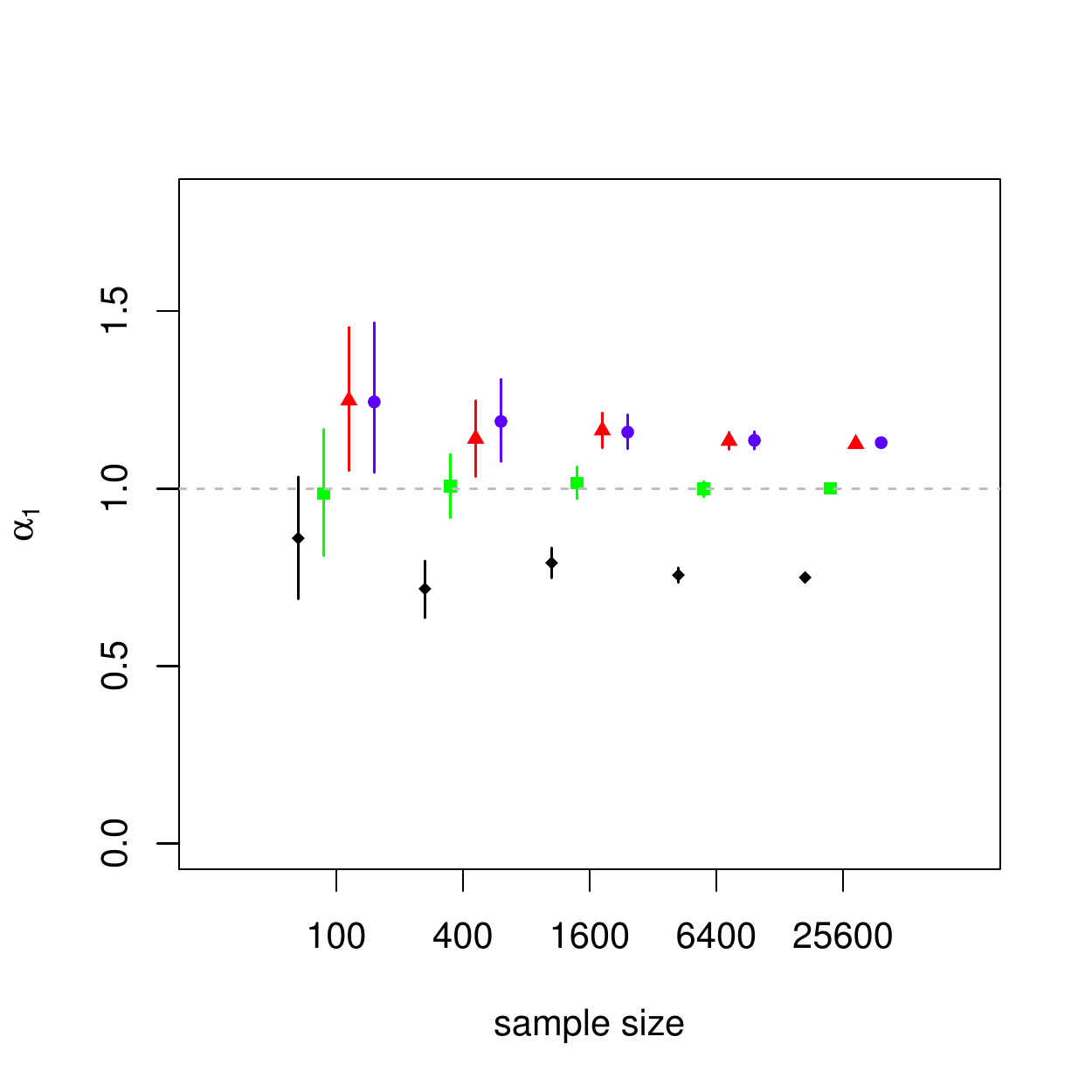,height=2.2in}}\quad
         \subfigure[$\alpha_0=1.2$, $w=10\%$]{\epsfig{figure=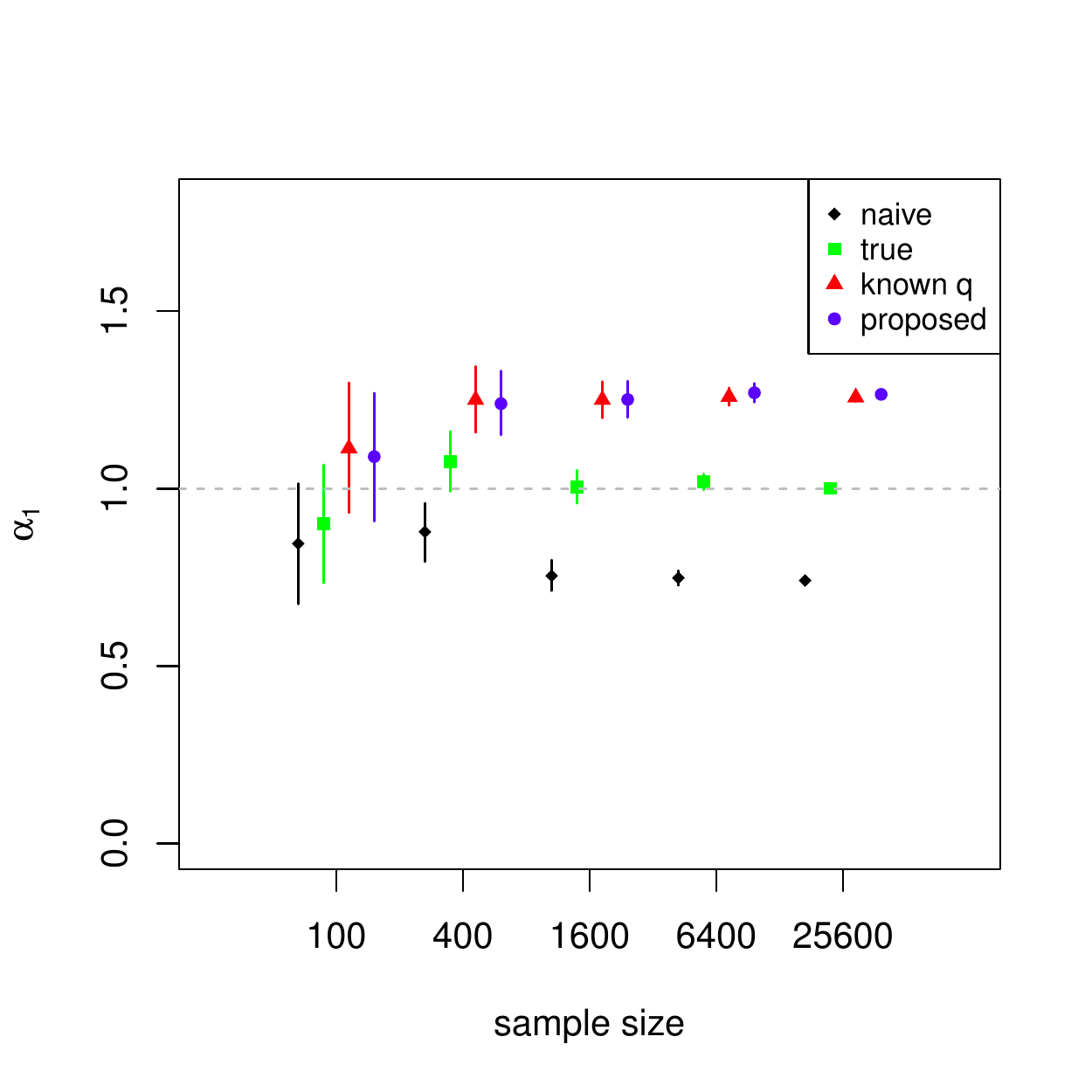,height=2.2in}}\quad
         \subfigure[$\alpha_0=0$, $w=5\%$]{\epsfig{figure=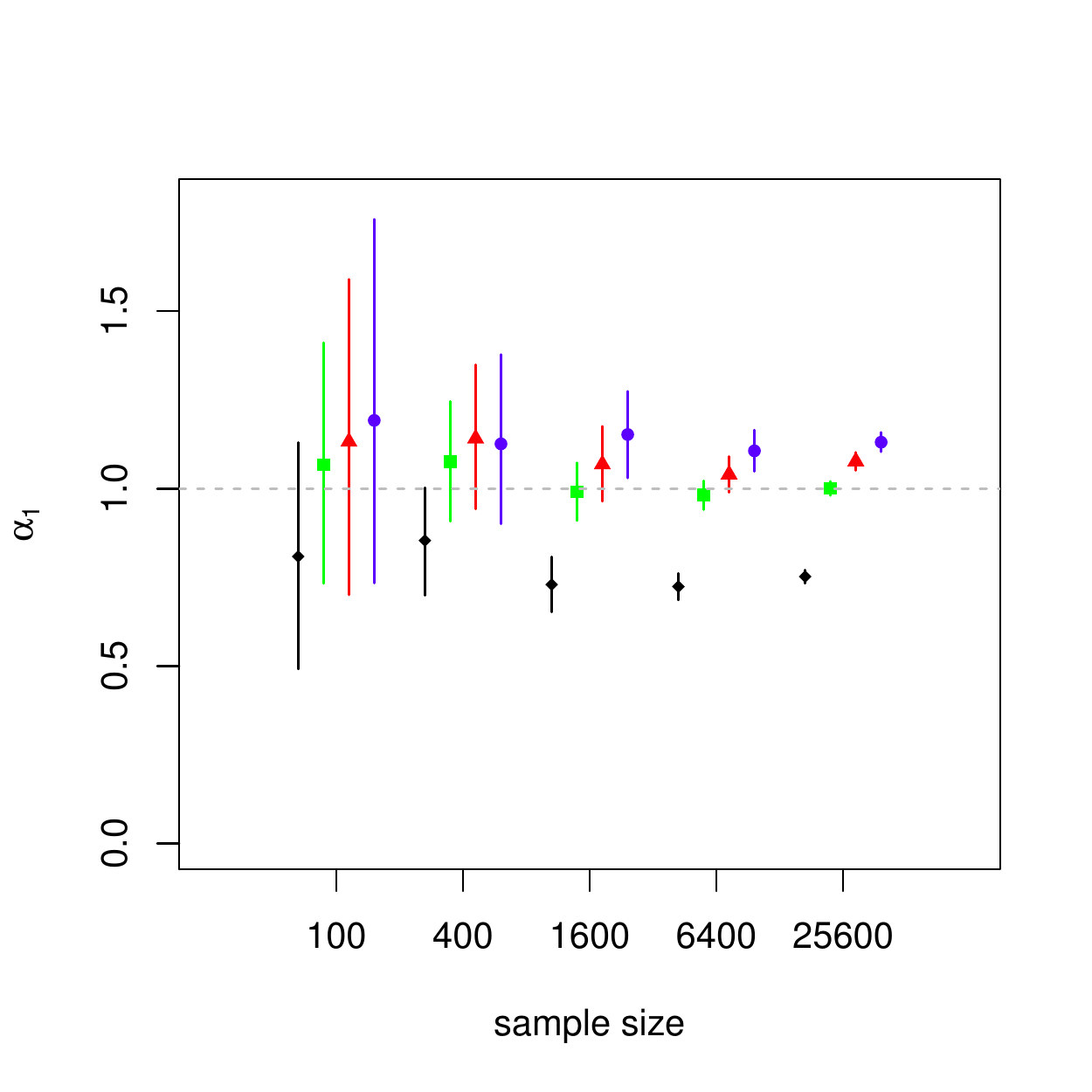,height=2.2in}}}     
     \caption{Equal-tailed credible intervals of $\alpha_1$ for ZIP responses with a binary $V$. The first panel corresponds to a variance to mean ratio of 1.19 for the $V=0$ case and 1.47 for the $V=1$ case. The second panel corresponds to a variance to mean ratio of 1.36 for the $V=0$ case and 2.01 for the $V=1$ case. The last panel corresponds to a smaller effect size, with a variance to mean ratio of 1.06 for the $V=0$ case and 1.15 for the $V=1$ case.}
   \label{poisson.sen1}
\end{figure}

Figure \ref{poisson.sen1} shows that the effect is over-estimated for the mixture models regardless whether the reclassification probabilities are known or not. It seems the two methods treat some or all of the additional zeros as coming from the $V=0$ group that has a smaller mean. The main observations are: (1) the bias increases with the severity of zero inflation; (2) knowing the reclassification offers more benefit in the case with a smaller effect size; (3) the reallocation of the additional zeros seems to have led to a small decrease in the variability of estimation, when comparing to the results in the earlier subsection. 

\section{Empirical analysis: Lung HIV study}\label{appl}
\subsection{Data}
The Longitudinal Studies of HIV-Associated Lung Infections and Complications (Lung HIV, \cite{crothers2011hiv}) was a collaborative multi-site study conducted between 2007 and 2013 by the National Heart, Lung, and Blood Institute (NHLBI). The study undertook data and specimen collection from eight HIV and pulmonary studies associated with NHLBI. The study's intent was to advance knowledge on HIV-related pulmonary diseases, motivated by the high incidence of serious pulmonary complications in HIV patients. 

The Lung HIV study included adults over the age of 18 who were diagnosed of HIV and self-reported smoking. There were a total of 904 participants included in the study who had a computed tomography (CT) performed at the baseline. In particular, the CT provided lung density measures including the relative area of emphysema (RA) below $-$910 Hounsfield units (HU, RA$-$910) and the 15th percentile density in HU (PD15) that have been regarded as effective assessments of the extent and severity of chronic obstructive pulmonary disease (COPD, see, e.g., \cite{soejima2000longitudinal,shaker2011rapid}). Self-reported traits on use of recreational drugs such as cocaine and heroin were collected on a voluntary basis. Use of illegal drugs, particularly cocaine, is a known risk factor for pulmonary disease in the general population (see, e.g., \cite{yakel1995pulmonary,alnas2010clinical,simonetti2014pulmonary,fiorelli2016spontaneous}). Baseline demographic information including age, gender, race, education, work status, family size and living arrangement are also included. 

After excluding records with missing values in the corresponding responses and covariates, the sample size ranges from 436 to 442 for the six distinct response variables. The numbers of participants who reported a positive status are 400 (male), 346 (exposure to smoking), 173 (white), 263 (cocaine use) for the binary traits. The age of the participants ranges from 21 to 75, and the number of cigarettes each participant smokes daily ranges from 0 to 45. The education variable is ordinal with 6 categories, with the value increasing with the level of education. We perform a logarithms transformation for the RA$-910$ and RA $-600$ to $-250$ variables, and all the variables seem to satisfy the normality assumption after the transformation. 

\subsection{Related literature}
Several recent papers have analyzed data from the Lung HIV study or its sub-studies, focused on understanding pulmonary complications related to HIV infection. Among them, \cite{drummond2015factors} and \cite{leader2016risk} appear to be the only two that analyzed the multi-center Lung HIV data. Both papers study risk factors such as age and CD4 cell count for lung complications in HIV-infected individuals. Additional papers use data from sub-studies of the broader Lung HIV study. Among them, papers such as \cite{sigel2014findings} and \cite{depp2016risk} study the association between HIV infection and pulmonary diffusing capacity and lung density measures, using data from the Examination of HIV-Associated Lung Emphysema and Veterans Aging Cohort Study that included both HIV-infected and non-infected participants. \cite{drummond2013effect} and \cite{lambert2015hiv} assess the association between HIV infection and lung function decline among injection drug users, using the AIDS Linked to the IntraVenous Experience cohort data on injection drug users. 

Of these previous papers, the most similar to our analysis is \cite{simonetti2014pulmonary}, which examines the association between pulmonary function and drug use in HIV-infected individuals, using data from the Multicenter AIDS Cohort Study and Women's Interagency HIV Study cohorts. They found no significant effect of recreational drug use on pulmonary function. There are two factors that may have limited the study's ability to detect a significant drug effect. First, the study had a total sample of only 184 patients, with 84 drug users. This sample size is small for the binary logistic regression model that was used. Second, intentional misrepresentation on drug use status was unaccounted for, although likely, due to legal/social desirability concerns.

Our approach addresses both of these concerns. By using the multi-center Lung HIV data we have a larger sample size, and our approach naturally handles the possibility of misclassification of the drug use variable.

\subsection{Analysis and results}
We are interested in assessing the effect of cocaine use in a continuous regression model on measures of lung density, after accounting for other risk factors and potential misclassification of the cocaine use variable. The response variables of interest include  lung density, PD15, mean density, RA below $-$910, RA below $-$856, and RA between $-600$ to $-250$ HU. Among the demographic and risk factors available in the study, we include age, gender, ethnicity, exposure to smoking (ExpSmk), the number of cigarettes smoked each day (NoCigs), education (Edu), and use of crack cocaine (CraCoc) as covariates. 

We perform an adjusted analysis using the proposed approach, and an unadjusted analysis using regular Bayesian normal model. We implement both the unadjusted and adjusted models in WinBUGS, through the R package \textbf{R2WinBUGS}. For all models, we run three chains with thinning of 6 and a burn-in of 15,000. Normal priors with mean 0 and variance 100 are used for the regression coefficients, except for the intercepts of the PD15 model, for which we set the prior variance to 10,000 due to its large scale. For the normal precision variable, we specify a vague gamma prior with parameters 1.5 and 2. For the reclassification probabilities, $q_{01}$ and $q_{10}$, we assume they do not depend on the other covariates. We use a beta prior with parameters 1 and 9, corresponding to a prior mean and standard deviation of 0.1. The prior puts more weights at the values near zero, so we are assuming the chance of misclassification is small, unless data suggest otherwise. As suggested by \cite{ferrari2008bayesian}, evidence from more informative studies may be incorporated through priors in a Bayesian model when such information is available. For all parameters, we take a posterior sample of size 5,000 with an effective size over 4,500.  

The prior and posterior densities for the reclassification probabilities are presented in Web Figure 25. the posterior mean and standard deviation for the reclassification probabilities ($q_{01}$, $q_{10}$) are [0.67(0.05), 0.22(0.04)] for the RA$-910$ model. The posterior standard deviation is much smaller than the prior standard deviation, indicating significant learning in the reclassification probabilities. The RA$-910$ model indicates a much larger probability of false negatives than that of false positives. This is expected due to the social undesirability and legal concern surrounding recreational drug use. The posterior mean of $q_{01}$ from the PD15 model is slightly larger (0.11) than the prior mean (0.1), indicating very weak learning of the reclassification probability when the effect size is small. For the density model, the posterior mean and standard deviation of the reclassification probabilities are the same as those of the priors. This is presumably due to the nonidentifiability of the model when the effect size of the cocaine use variable is zero. 

\begin{figure}[htp]
   \centering
   \mbox{\subfigure[RA$-$910]{\epsfig{figure=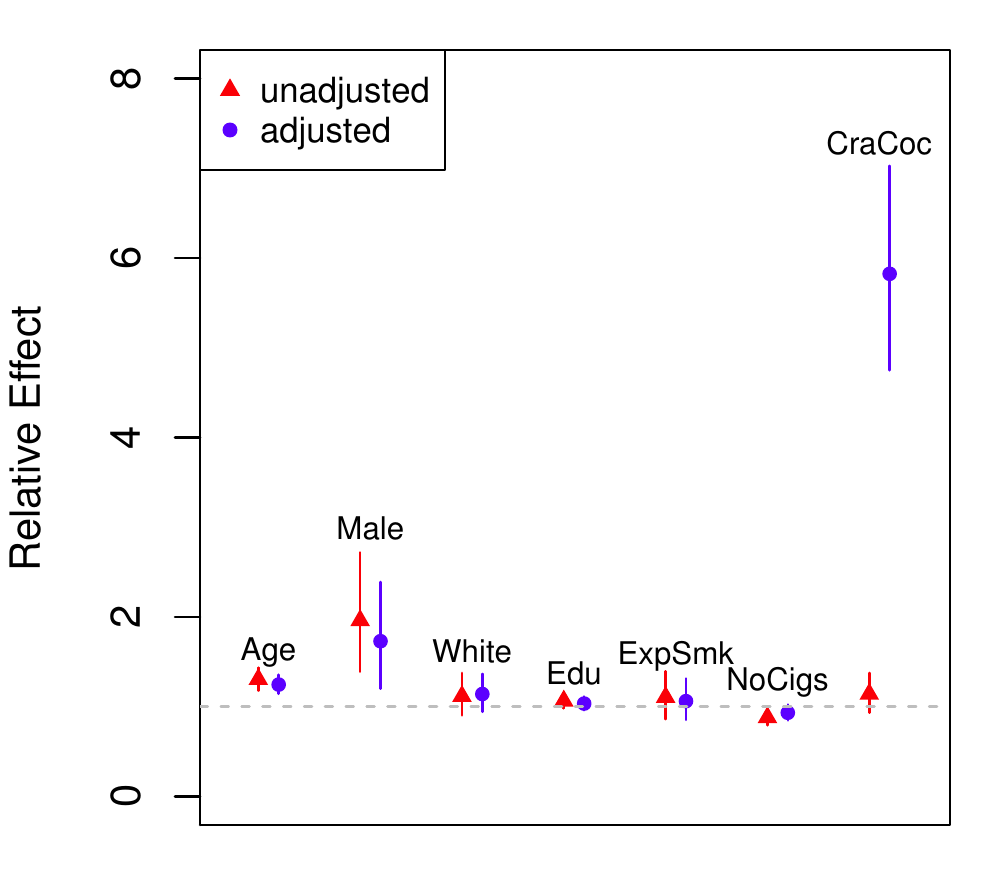,height=1.9in}}\quad
         \subfigure[PD15]{\epsfig{figure=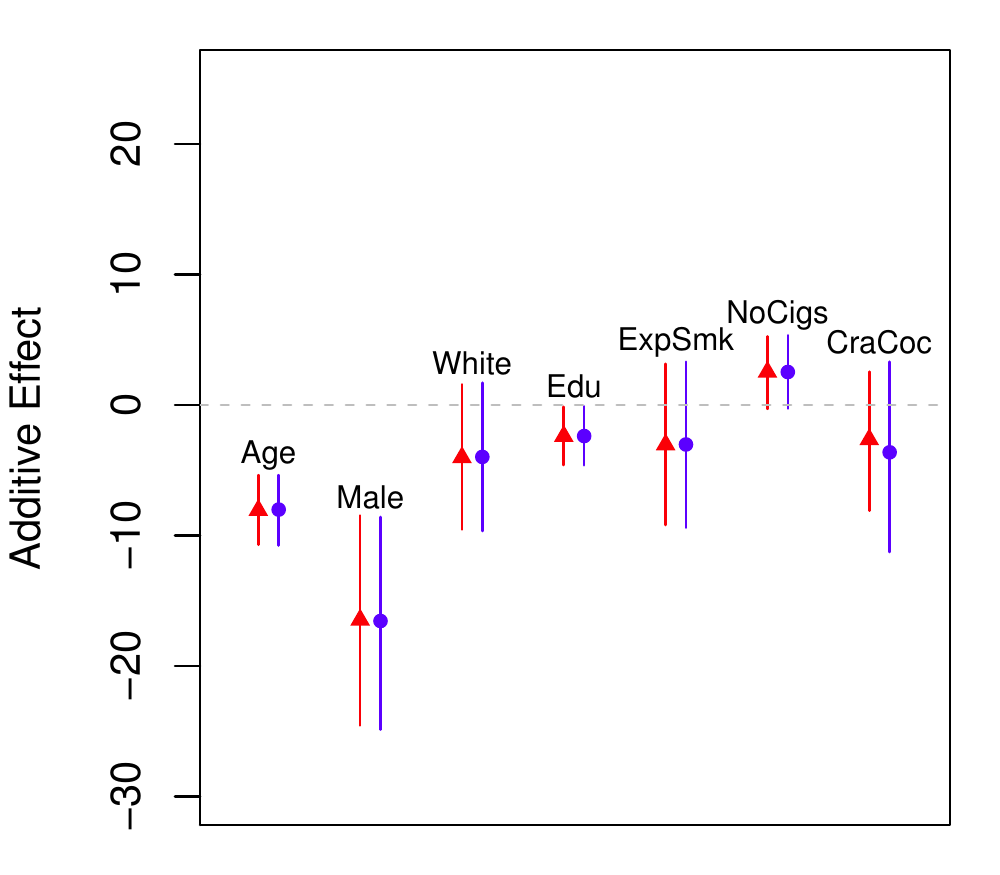,height=1.9in}}\quad
         \subfigure[Density]{\epsfig{figure=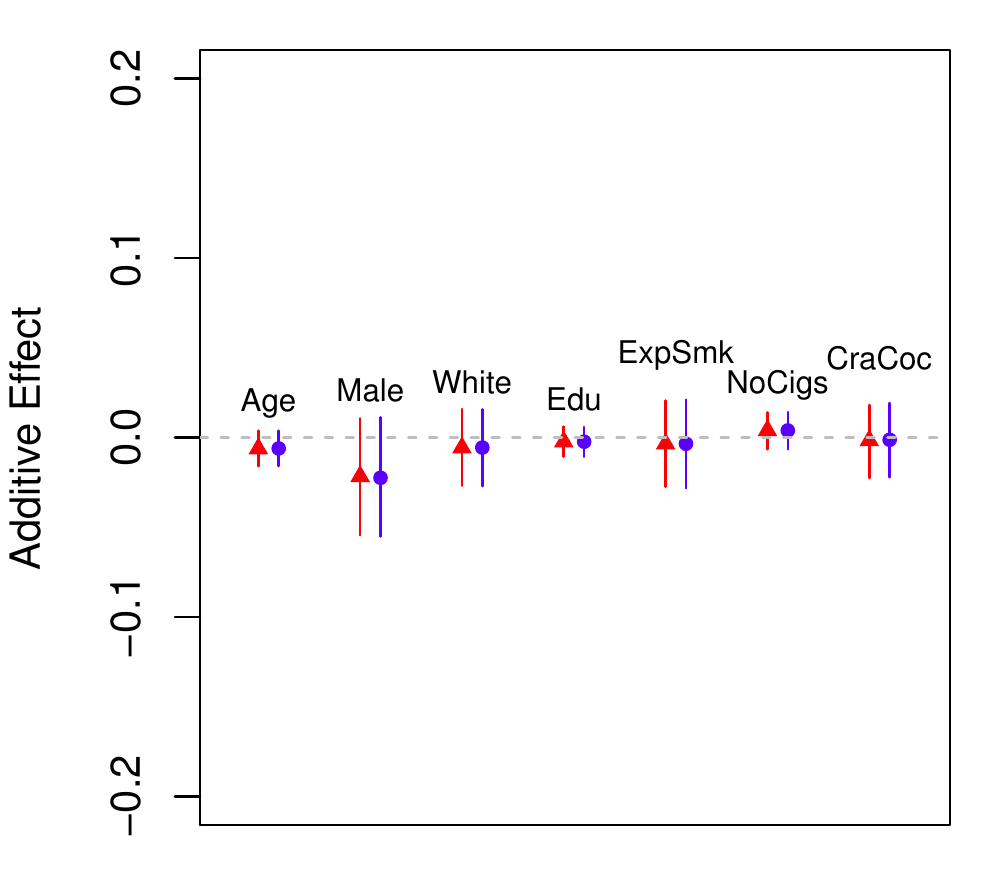,height=1.9in}}}
     \caption{Equal-tailed credible intervals of covariate effects on the lung density measures including the RA$-$910, PD15, and lung density. A log transformation is performed on RA$-910$ in order to satisfy the normality assumption. Hence, the exponential of the coefficient corresponds to the relative effect on the median of RA$-910$. The age effect corresponds to every 10-year increase in the age, and the cigarette effect corresponds to every 10-cigarette increase per day.}
   \label{hiv}
\end{figure}

The 95\% equal-tailed credible intervals for the covariate effects are given in Figure \ref{hiv} and Web Figure 26. For the two variables of RA$-$910 and RA $-600$ to $-250$, we fit the models on the logarithm scale, so we present the exponential of the slope that represents the multiplicative effect on the median of these variables. For the other variables, we present the regression slope that can be interpreted as the additive effect on the mean of the variables.

In Figure \ref{hiv}, we observe that the directions of the effect for all covariates match those found in the previous literature. Of particular note is that the number of cigarettes a participant smoked has an estimated effect that is opposite to that of the age and cocaine use variables. This is consistent with the findings from \cite{shaker2011rapid}: inflammation from smoking may mask the presence of emphysema on CT. After adjusting for misclassification in the self-reported cocaine use status, we observe a significant effect of cocaine use on worsening lung complications measure by the RA below $-$910 HU. The estimated effect size for the model is around 2.6, confirming it is a case where the proposed model is efficient. Since the relative effect is the exponential of the regression coefficient, the adjustment results in a large impact on the estimated relative effect. For the other two response variables, PD15 and lung density, the adjustment does not result in a significant difference in the estimated effects. 

In summary, our analysis suggests that the high incidence of cocaine use in HIV-infected individuals is probably a major contributing factor for the severe or fatal pulmonary complications that have been observed in the target population. 

\subsection{Sensitivity analysis}
Here, we perform a sensitivity analysis on how the results concerning RA$-$910 can be affected when the reclassification probabilities are permitted to vary as a function of other covariates.

Specifically, we model the reclassification probabilities $q_{01}$ and $q_{10}$ using the logit model given in (\ref{mix.eq3}). In the binary logit model, we include covariates race, education and age, which are plausibly associated with cocaine use or misclassification status. We investigate six logit models including either one, two, or all three covariates.

\begin{figure}[htp]
   \centering
   \epsfig{figure=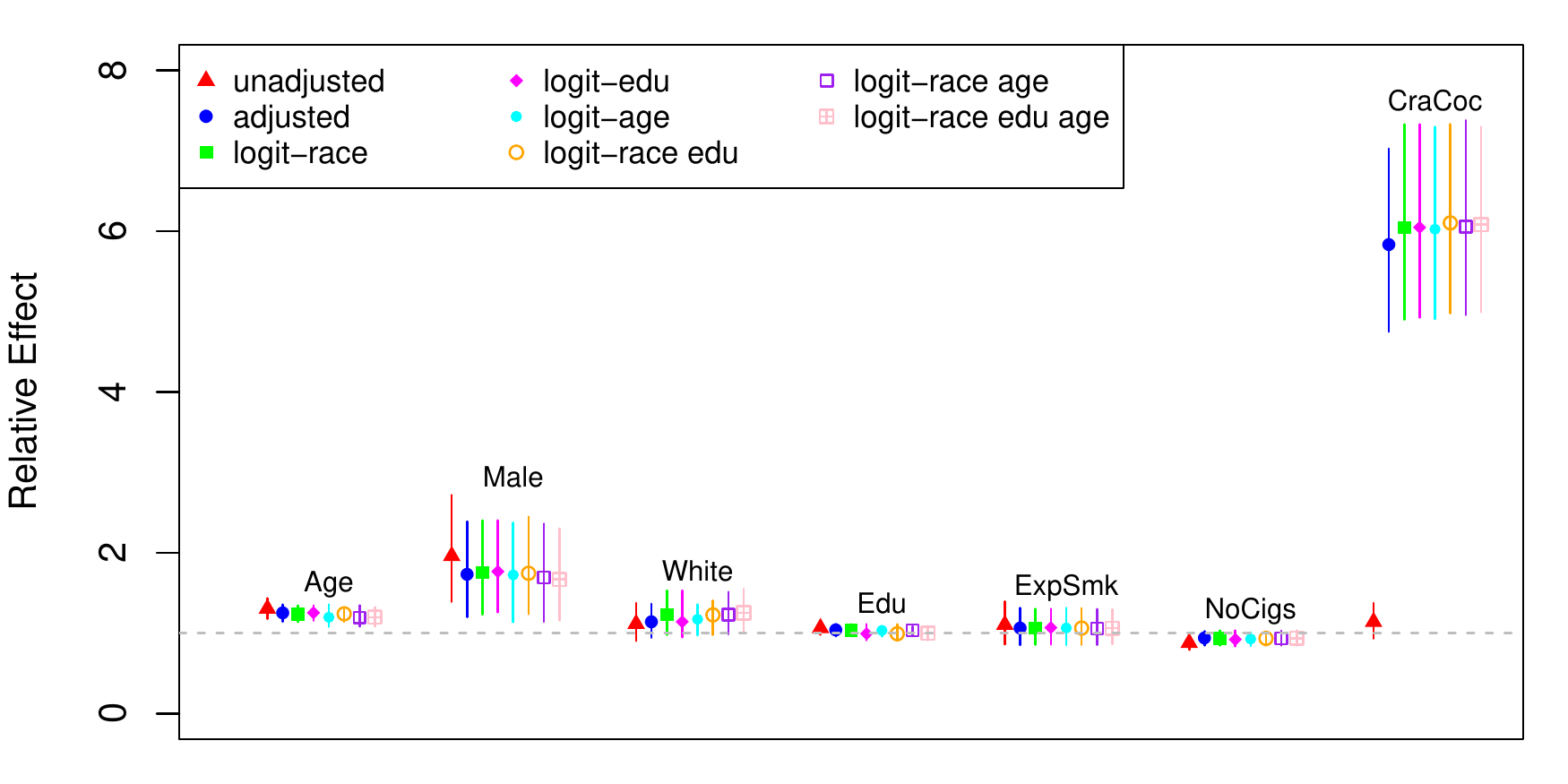,height=2.8in}
     \caption{Equal-tailed credible intervals of covariate effects on RA$-$910, for the sensitivity analysis with logit models on reclassification probabilities $q_{01}$ and $q_{10}$. }
   \label{hiv2}
\end{figure}

MCMC and prior settings are similar to those in the previous subsection. In particular, for the regression coefficients of the logit models, we select the priors in an effort to match the unconditional prior mean and standard deviation for reclassification probabilities, those of the $beta(1,\,9)$ prior we assumed in the previous subsection. For the single-covariate logit models, we standardize the corresponding covariates before using a $N(-2.5,\,0.04)$ prior for the intercept and a standard normal prior for the slope. For the multi-covariate logit models, the priors are selected in a similar manner. Figure \ref{hiv2} shows the 95\% equal-tailed credible intervals of the covariate effects on RA$-$910 for the models we have considered. The prior and posterior distributions of the regression coefficients for the logit models are presented in Web Figure 27.

Figure \ref{hiv2} shows that regression on reclassification probabilities leads to a mild increase in the estimated cocaine effect, with the cost of a mild increase or decrease in those of the corresponding covariates. From Web Figure 27, none of the covariate effects seems to differ significantly from zero. Overall, the basic findings appear to be robust to the modeling assumptions governing the reclassification probabilities.

\section{Conclusions}\label{concl}
In this paper, we formulate regression models with a misclassified categorical covariate as mixture regression models. The finite mixture representation enables us to perform valid statistical inference on all parameters in the model, including the regression coefficients and reclassification probabilities, without requiring extra sources of information on the misclassification probabilities or the true covariate value. By evaluating the efficiency loss caused by the misclassification, we showed that the efficiency loss is dominated by the effect size, the severity of misclassification, and the distribution for the categorical covariate of concern. For cases where the effect size is large, not observing the true covariate value and not knowing the reclassification probabilities contributes little to loss of efficiency. Therefore, conducting inference using the mixture represetnation allows for reliable inference on model parameters even when there is no information or validation data available on the misclassification or the true status. Further, furnishing side information on the reclassification probabilities via an informative prior protects against the possibility that the side information is wrong.
Applying the methodology to the Lung HIV study data, we confirmed that the adjustment of misclassification in the self-reported cocaine use results in an estimated effect that is substantially larger than that indicated by the unadjusted analysis. Our analysis indicates a significant effect of crack cocaine use in worsening lung complications measured by the relative area of emphysema, after adjusting for other known risk factors. 

\section*{Acknowledgments}
This manuscript was prepared using Lung HIV Research Materials obtained from the NHLBI Biologic Specimen and Data Repository Information Coordinating Center and does not necessarily reflect the opinions or views of the Lung HIV and the NHLBI.

\bibliographystyle{biom}


\bibliography{biblio}

\end{document}